\documentclass[a4paper,11pt]{article}

\usepackage[utf8]{inputenc}

\usepackage[T1]{fontenc}
\usepackage{fouriernc}

\usepackage{amsmath}
\usepackage{amssymb}
\usepackage{bbding}
\usepackage{color}
\usepackage{enumitem}
\usepackage{graphicx}
\usepackage{nicefrac}
\usepackage{ntheorem}

\usepackage{tikz}
\newcommand*\circled[1]{\tikz[baseline=-3.8pt]{\node[shape=circle,draw,inner sep=0.5pt] (char) {\scriptsize #1};}}

\usepackage[colorlinks=true,linkcolor=rot,citecolor=rot,urlcolor=rot]{hyperref}

\usepackage[mathlines]{lineno}

\definecolor{grau}{RGB}{128,128,128}
\definecolor{rot}{RGB}{179,27,27}

\newcommand{\schummeln}[1]{\hspace{0.#1mm}}
\newcommand{\wegschummeln}[1]{\hspace{-0.#1mm}}

\newcommand{\llangle}{\left\langle\hspace{-1.5pt}\left\langle}
\newcommand{\rrangle}{\right\rangle\hspace{-1.5pt}\right\rangle}
\newcommand{\lllangle}{\left\langle\hspace{-2.25pt}\left\langle}
\newcommand{\rrrangle}{\right\rangle\hspace{-2.25pt}\right\rangle}

\newcommand{\pathbrackets}[1]{[\hspace{-2pt}[#1]\hspace{-2pt}]}

\DeclareMathOperator{\dist}{d}
\DeclareMathOperator{\interior}{int}
\DeclareMathOperator{\isom}{Isom}
\DeclareMathOperator{\wt}{wt}

\newtheorem{definition}{Definition}[section]
\newtheorem{corollary}[definition]{Corollary}
\newtheorem{lemma}[definition]{Lemma}
\newtheorem{proposition}[definition]{Proposition}
\newtheorem{remark}[definition]{Remark}

\newtheorem*{assumption}{Standing Assumption}
\newtheorem{theorem}{Theorem}

\newcounter{recallenumi}

\title{The Tits alternative for\\non-spherical triangles of groups}

\author{
Johannes Cuno\thanks{The first author acknowledges the support of the Austrian Science Fund (FWF), project W1230-N13, and the support of the Cornell Graduate School.}\\
Technische Universit\"at Graz\\
\texttt{cuno@math.tugraz.at}\\[8pt]
J\"org Lehnert\\
Max Planck Institute for\\ Mathematics in the Sciences\\
\texttt{lehnert@mis.mpg.de}
}

\date{}

\begin{document}

% Zeilenangabe zum Korrekturlesen:
% \linenumbers

% ==================================================
% TITELSEITE, ABSTRACT UND INHALTSVERZEICHNIS
% ==================================================

\maketitle

\begin{abstract}
\noindent Triangles of groups have been introduced by Gersten and Stallings. They are, roughly speaking, a generalisation of the amalgamated free product of two groups and occur in the framework of Corson diagrams. First, we prove an intersection theorem for Corson diagrams. Then, we focus on triangles of groups. It has been shown by Howie and Kopteva that the colimit of a hyperbolic triangle of groups contains a non-abelian free subgroup. We give two natural conditions, each of which ensures that the colimit of a non-spherical triangle of groups either contains a non-abelian free subgroup or is virtually solvable. \smallskip

% ==================================================
% STICHWÖRTER
% ==================================================

 \smallskip

% ==================================================
% AMS-KLASSIFIKATION
% ==================================================
% 20E06 - Free products, free products with
%         amalgamation, Higman-Neumann-Neumann
%         extensions, and generalizations
% 20F05 - Generators, relations, and presentations
% 20F65 - Geometric group theory
% 20E05 - Free nonabelian groups
% 20E07 - Subgroup theorems; subgroup growth
% ==================================================

\noindent \textbf{MSC 2010 classes:} 20E06 (primary), 20F05, 20F65, 20E05, 20E07.
\end{abstract}

% ==================================================
% §1 INTRODUCTION
% ==================================================
%
\section{Introduction} \label{sec:introduction}
Given a commutative diagram of groups and injective homomorphisms, we may construct its colimit (in the category of groups). The colimit, or, more precisely, the colimit group, can be obtained by taking the free product of the groups and identifying the factors according to the homomorphisms. A good example is the amalgamated free product $X\ast_{A} Y$, which is the colimit group of the diagram $X\leftarrow A\rightarrow Y$.

\begin{figure}
\begin{center}
\input{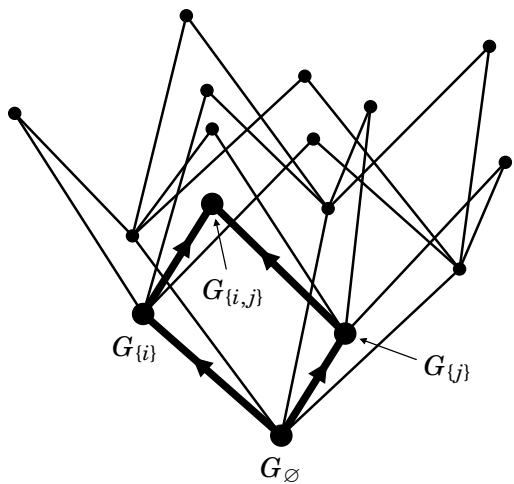}
\end{center}
\caption{A Corson diagram based on a set $I$ with $|I|=5$. For simplicity, the homomorphisms $\varphi_{\emptyset J}:G_{\emptyset}\rightarrow G_{J}$ with $J\subseteq I$ and $|J|=2$ have been omitted.}
\label{fig:corson}
\end{figure}

A Corson diagram is based on a set $I$. For every subset $J\subseteq I$ with $|J|\leq 2$ there is a group $G_{J}$ and for every two subsets $J_{1}\subset J_{2}\subseteq I$ with $|J_{2}|\leq 2$ there is a homomorphism $\varphi_{J_{1}J_{2}}:G_{J_{1}}\rightarrow G_{J_{2}}$, see Figure~\ref{fig:corson}.
Notice that both Artin groups and Coxeter groups have a natural interpretation as colimit groups of Corson diagrams. A triangle of groups is nothing but a Corson diagram based on a set $I$ with $|I|=3$. Gersten and Stallings introduced the notion of curvature and proved that for non-spherical triangles of groups the natural homomorphisms $\nu_{J}$ from the groups $G_{J}$ to the colimit group $\mathcal{G}$ are injective, see \cite{Stallings1991}. A similar result holds for non-spherical Corson diagrams, see \cite{Corson1996}. While these results can be proved by nice arguments based on Euler's formula for planar graphs, spherical Corson diagrams are much harder to investigate, see e.\,g.~\cite{Chermak1995} and \cite{Allcock2012}.

In Section~\ref{sec:preliminaries}, we introduce some basic notions. Then, in Section~\ref{sec:intersection}, we give an example of a spherical triangle of groups showing that, even though the natural homomorphisms $\nu_{J}:G_{J}\rightarrow\mathcal{G}$ are injective, the intersections of their images may be larger than the amalgamated subgroups. On the other hand, this cannot happen for non-spherical triangles of groups and, more generally, non-spherical Corson diagrams, see Theorem~\ref{thm:intersection}.
\begin{remark}
The absence of large intersections shall not be confused with the developability of complexes of groups, which is implied by the injectivity of the natural homomorphisms $\nu_{J}:G_{J}\rightarrow\mathcal{G}$, see \cite[Corollary III.C.2.15]{BridsonHaefliger1999}.
\end{remark}
Howie and Kopteva showed that, under mild assumptions, the colimit group of a hyperbolic triangle of groups has a non-abelian free subgroup, see \cite{HowieKopteva2006}. In Section~\ref{sec:billiards}, we focus on the Euclidean case and discuss the following version of the Tits alternative: ``A class $\mathfrak{C}$ of groups satisfies the Tits alternative if each $G\in\mathfrak{C}$ is either large, i.\,e.~has a non-abelian free subgroup, or small, i.\,e.~is virtually solvable.'' The Tits alternative is named after Jacques Tits, who proved in 1972 that a finitely generated linear group is either large or small, see \cite[Corollary~1]{Tits1972}. Since then, the Tits alternative has been proved for many other classes of groups. For a list of results and open problems we refer to \cite{SageevWise2005}.

As indicated above, we are interested in Euclidean triangles of groups. In the case that none of the Gersten-Stallings angles is 0, we may follow Bridson's construction of a simplicial complex $\mathcal{X}$, see \cite{Bridson1991}, and use billiards on a triangle $\Delta\subseteq\mathbb{E}^{2}$ to obtain geodesics in the geometric realisation $|\mathcal{X}|$. These geodesics allow us to prove that, as soon as the simplicial complex $\mathcal{X}$ branches, the colimit group has a non-abelian free subgroup, see Theorems~\ref{thm:billiard} and \ref{thm:freesubgroup}.

The remaining cases can be analysed with quotients and amalgamated free products. In the end, we generalise the result by Howie and Kopteva mentioned above and prove that the Tits alternative holds for the class of colimit groups of non-spherical triangles of groups with the property that:
\begin{itemize}
\item none of the Gersten-Stallings angles is 0 and the group $G_{\emptyset}$ is finitely generated and either large or small, see Theorem~\ref{thm:titsone}, or
\item every group $G_{J}$ with $J\subseteq I$ and $|J|\leq 2$ is finitely generated and either large or small, see Theorem~\ref{thm:titstwo}.
\end{itemize}
\subsection*{Acknowledgement}
This paper originates in a Diplomarbeit, the equivalent of a master's thesis, under supervision of Robert Bieri at Goethe-Universit\"at Frankfurt am Main. In the light of this beginning, we would like to thank Robert Bieri for his enthusiasm, advice, and patience.
%
% ==================================================
% §2 PRELIMINARIES
% ==================================================
%
\section{Preliminaries} \label{sec:preliminaries}
\subsection{Corson diagrams and their colimits} \label{sub:playground}
Let $I$ be an arbitrary set. Assume given for every subset $J\subseteq I$ with $|J|\leq 2$ a group $G_{J}$ and for every two subsets $J_{1}\subset J_{2}\subseteq I$ with $|J_{2}|\leq 2$ an injective homomorphism $\varphi_{J_{1}J_{2}}:G_{J_{1}}\rightarrow G_{J_{2}}$. Moreover, assume the resulting diagram to be \emph{commutative}, i.\,e.~for every sequence of subsets $\emptyset=J_{1}\subset J_{2}\subset J_{3}\subseteq I$ with $|J_{3}|=2$ the equation $\varphi_{J_{1}J_{3}}=\varphi_{J_{2}J_{3}}\circ\varphi_{J_{1}J_{2}}$ holds.

Since these diagrams have been introduced by Corson in \cite{Corson1996}, we refer to them as \emph{Corson diagrams}. In the case $|I|=3$, Corson diagrams are known as \emph{triangles of groups}. Whenever we consider a triangle of groups, we may assume w.\,l.\,o.\,g.~that $I=\{1,2,3\}$.

Given a Corson diagram, we will mostly be interested in its \emph{colimit group}. The colimit group can be obtained by taking the free product of the groups $G_{J}$ and identifying the factors according to the homomorphisms. Let us make this construction a little more precise! Think of the groups $G_{J}$ to be given by presentations $G_{J}\cong\langle G_{J}:R_{J}\rangle$, where $R_{J}$ is the set of all words over the group elements and their formal inverses that represent the identity. Then, the colimit group $\mathcal{G}$ is given by the following presentation:
\begin{equation}
\label{def:colimit}
\tag{$\ast$}
\mathcal{G}\cong\left\langle\bigsqcup_{\genfrac{}{}{0pt}{}{J\subseteq I}{|J|\leq 2}}G_{J}\,:\,\bigsqcup_{\genfrac{}{}{0pt}{}{J\subseteq I}{|J|\leq 2}}R_{J},\bigsqcup_{\genfrac{}{}{0pt}{}{J_{1}\subset J_{2}\subseteq I}{|J_{2}|\leq 2}}\left\{g=\varphi_{J_{1}J_{2}}(g)\,:\,g\in G_{J_{1}}\right\}\right\rangle
\end{equation}
This presentation, though not very economic, turns out to be suitable for our purposes. For every subset $J\subseteq I$ with $|J|\leq 2$ we may consider the natural homomorphism $\nu_{J}:G_{J}\rightarrow\mathcal{G}$ given by $g\mapsto g$. The colimit group, equipped with these homomorphisms, is called the \emph{colimit}. For further reading about it we refer to \cite[1.1]{Tits1986} and \cite[Chapter 3]{AdamekHerrlichStrecker1990}.
\subsection{Curvature of Corson diagrams} \label{sub:curvature}
The homomorphisms $\nu_{J}:G_{J}\rightarrow\mathcal{G}$ do not need to be injective. An example of a triangle of groups in which they are not has been given by Gersten and Stallings in \cite[1.4]{Stallings1991}. On the other hand, it turned out that for non-spherical triangles of groups and, more generally, for non-spherical Corson diagrams they are. Let us therefore introduce the notion of curvature! For every two distinct $i,j\in I$ the homomorphisms $\varphi_{\{i\}\{i,j\}}$ and $\varphi_{\{j\}\{i,j\}}$ uniquely determine a homomorphism $\alpha:G_{\{i\}}\ast_{G_{\varnothing}}G_{\{j\}}\rightarrow G_{\{i,j\}}$. If $\alpha$ is not injective, let $\hat{m}$ denote the minimal length of a non-trivial element in its kernel (in the usual length function on the amalgamated free product). Notice that the homomorphisms $\varphi_{\{i\}\{i,j\}}$ and $\varphi_{\{j\}\{i,j\}}$ are injective, whence the minimal length $\hat{m}$ must be even. The \emph{Gersten-Stallings angle} $\sphericalangle_{\{i,j\}}$ is now defined by:
\[
\sphericalangle_{\{i,j\}}=\left\{\begin{array}{cl} \nicefrac{2\pi\,}{\hat{m}} & \text{if $\alpha$ is not injective} \\  0 & \text{if $\alpha$ is injective} \end{array}\right.
\]
Three pairwise distinct elements $i,j,k\in I$ are called a \emph{spherical triple} if the sum $\sphericalangle_{\{i,j\}}+\sphericalangle_{\{i,k\}}+\sphericalangle_{\{j,k\}}$ is strictly larger than $\pi$. The Corson diagram is called \emph{spherical} if it has a spherical triple, and \emph{non-spherical} otherwise.
\begin{remark}
Consider a non-spherical triangle of groups. Since $I=\{1,2,3\}$, there is only one set of three pairwise distinct elements. Depending on whether the sum $\sphericalangle_{\{1,2\}}+\sphericalangle_{\{1,3\}}+\sphericalangle_{\{2,3\}}$ is strictly smaller than $\pi$ or equal to $\pi$, the \mbox{triangle} of groups is called hyperbolic or Euclidean. This distinction is of relevance in Section~\ref{sec:billiards}.
\end{remark}
\subsection{Embedding theorems} \label{sub:corson}
We are now able to state the theorem about non-spherical triangles of groups that has been mentioned above.
\begin{theorem}[Gersten-Stallings]
\label{thm:gs}
For every non-spherical triangle of groups the natural homomorphisms $\nu_{J}:G_{J}\rightarrow\mathcal{G}$ are injective.\hfill $\Box$
\end{theorem}
Theorem~\ref{thm:gs} has been proved in \cite{Stallings1991} and generalised to non-spherical Corson diagrams in \cite{Corson1996}. Even more has been shown in \cite{Corson1996}: Not only the groups $G_{J}$ but also the colimit groups of subdiagrams naturally embed into $\mathcal{G}$. Let us clarify! Given a Corson diagram and a subset $K\subseteq I$ we may restrict our focus to the subdiagram spanned by the groups $G_{J}$ with $J\subseteq K$ and $|J|\leq 2$, see the bold vertices and arrows in Figure~\ref{fig:corson} for an example. The colimit group of such a subdiagram can be obtained by modifying (\ref{def:colimit}) as follows:
\[
\mathcal{G}|_{K}\cong\left\langle\bigsqcup_{\genfrac{}{}{0pt}{}{J\subseteq K}{|J|\leq 2}}G_{J}\,:\,\bigsqcup_{\genfrac{}{}{0pt}{}{J\subseteq K}{|J|\leq 2}}R_{J},\bigsqcup_{\genfrac{}{}{0pt}{}{J_{1}\subset J_{2}\subseteq K}{|J_{2}|\leq 2}}\left\{g=\varphi_{J_{1}J_{2}}(g)\,:\,g\in G_{J_{1}}\right\}\right\rangle
\]
Analogously to $\nu_{J}:G_{J}\rightarrow\mathcal{G}$ introduced in Section~\ref{sub:playground}, we may now consider the natural homomorphisms $\tilde\nu_{K}:\mathcal{G}|_{K}\rightarrow\mathcal{G}$ given by $g\mapsto g$.
\begin{theorem}[Corson]
\label{thm:corson}
For every non-spherical Corson diagram the natural homomorphisms $\tilde\nu_{K}:\mathcal{G}|_{K}\rightarrow\mathcal{G}$ are injective.\hfill $\Box$
\end{theorem}
\begin{remark}
\label{rem:implication}
It is easy to verify that for every subset $K\subseteq I$ with $|K|\leq 2$ there is an isomorphism $\mu_{K}:G_{K}\rightarrow\mathcal{G}|_{K}$ given by $g\mapsto g$ and, therefore, the injectivity of $\tilde\nu_{K}:\mathcal{G}|_{K}\rightarrow\mathcal{G}$ implies the injectivity of $\nu_{K}=\tilde\nu_{K}\circ\mu_{K}:G_{K}\rightarrow\mathcal{G}$. In particular, Theorem~\ref{thm:corson} implies Theorem~\ref{thm:gs}.
\end{remark}
\begin{remark}
\label{rem:shorthand}
We make the following convention: Whenever we know that the homomorphisms $\tilde\nu_{K}:\mathcal{G}|_{K}\rightarrow\mathcal{G}$ are injective, e.\,g.~in case of a non-spherical Corson diagram, we do not need to mention them any more and may tacitly interpret $\mathcal{G}|_{K}$ as a subgroup of $\mathcal{G}$.
\end{remark}
\begin{remark}
\label{rem:shorthandinterpretation}
In this situation, the symbol $\mathcal{G}|_{K}$ refers to the subgroup of $\mathcal{G}$ that is generated by the elements of the groups $G_{J}$ with $J\subseteq K$ and $|J|\leq 2$. This way, we can easily observe that $K_{1}\subseteq K_{2}$ implies $\mathcal{G}|_{K_{1}}\subseteq\mathcal{G}|_{K_{2}}$.
\end{remark}
\subsection{Standing assumption on the Gersten-Stallings angles} \label{sub:assumption}
We will have to make one more assumption, which has already been indicated by Gersten and Stallings in \cite[p.\,493, ll.\,4--6]{Stallings1991} and Corson in \cite[p.\,567, l.\,15]{Corson1996}, even though Theorems~\ref{thm:gs} and \ref{thm:corson} hold without it.
\begin{assumption}
\label{ass:angles}
We shall always assume, without stating explicitly, that none of the Gersten-Stallings angles is equal to $\pi$.
\end{assumption}
\begin{remark}
\label{rem:anglesinterpretationpart1}
This assumption is equivalent to the property that for every two distinct $i,j\in I$ the equation $\varphi_{\{i\}\{i,j\}}\left(G_{\{i\}}\right)\cap\varphi_{\{j\}\{i,j\}}\left(G_{\{j\}}\right)=\varphi_{\emptyset\{i,j\}}\left(G_{\emptyset}\right)$ holds.
\end{remark}
%
% ==================================================
% §3 INTERSECTION THEOREM
% ==================================================
%
\section{Intersection theorem} \label{sec:intersection}
Assume given a Corson diagram with the property that the homomorphisms $\tilde\nu_{K}:\mathcal{G}|_{K}\rightarrow\mathcal{G}$ are injective. One question we are interested in is whether two subgroups $\mathcal{G}|_{K_{1}}$ and $\mathcal{G}|_{K_{2}}$ intersect only along the obvious subgroup $\mathcal{G}|_{K_{1}\cap K_{2}}$ or along some larger subgroup of $\mathcal{G}$. In Section~\ref{sub:example}, we give an example of a spherical Corson diagram in which the homomorphisms $\tilde\nu_{K}:\mathcal{G}|_{K}\rightarrow\mathcal{G}$ are injective but there are $K_{1},K_{2}\subseteq I$ such that $\mathcal{G}|_{K_{1}}\cap\mathcal{G}|_{K_{2}}\neq\mathcal{G}|_{K_{1}\cap K_{2}}$. Then, we recall the notion of \emph{disc pictures} and use it to prove an \emph{intersection theorem} showing that this can only happen in the spherical realm.
\subsection{Example} \label{sub:example}
Let us consider the Corson diagram based on the following data: $I=\{1,2,3\}$, $G_{\varnothing}\cong\{1\}$, $G_{\{1\}}\cong\left\langle a\,:\,-\right\rangle$, $G_{\{2\}}\cong\left\langle b\,:\,-\right\rangle$, $G_{\{3\}}\cong\left\langle c\,:\,-\right\rangle$, $G_{\{1,2\}}\cong\left\langle a,b\,:\,b^{-1}ab=a^{2}\right\rangle$, $G_{\{1,3\}}\cong\left\langle a,c\,:\,c^{-1}ac=a^{2}\right\rangle$, $G_{\{2,3\}}\cong\left\langle b,c\,:\,bc=cb\right\rangle$. Here, the homomorphisms $\varphi_{J_{1}J_{2}}:G_{J_{1}}\rightarrow G_{J_{2}}$ are implicitly given by $a\mapsto a$, $b\mapsto b$, $c\mapsto c$. Since $G_{\emptyset}$ is trivial, the resulting diagram is commutative. Britton's Lemma \cite[Principal Lemma]{Britton1963} shows that the homomorphisms $\varphi_{J_{1}J_{2}}:G_{J_{1}}\rightarrow G_{J_{2}}$ are injective and the Gersten-Stallings angles amount to $\nicefrac{\pi\,}{2}$ each. So, it is a spherical Corson diagram in the sense of Sections~\ref{sub:playground}, \ref{sub:curvature}, and \ref{sub:assumption}. 
\begin{proposition}
%\label{prn:examplepart1}
The natural homomorphisms $\tilde\nu_{K}:\mathcal{G}|_{K}\rightarrow\mathcal{G}$ are injective.
\end{proposition}
\textbf{Proof.} Since the homomorphism $\tilde\nu_{\{1,2,3\}}:\mathcal{G}|_{\{1,2,3\}}\rightarrow\mathcal{G}$ is obviously injective, it suffices to verify the injectivity of the homomorphisms $\tilde\nu_{K}:\mathcal{G}|_{K}\rightarrow\mathcal{G}$ with $K\subseteq\{1,2,3\}$ and $|K|\leq 2$. But then, we already know from Remark~\ref{rem:implication} that there are isomorphisms $\mu_{K}:G_{K}\rightarrow\mathcal{G}|_{K}$ with $\nu_{K}=\tilde\nu_{K}\circ\mu_{K}$. So, it even suffices to verify the injectivity of the homomorphisms $\nu_{K}:G_{K}\rightarrow\mathcal{G}$.

Recall the presentation (\ref{def:colimit}) of the colimit group $\mathcal{G}$ and notice that, in our situation, it can be simplified by deleting superficial generators and relators so that we finally obtain the following presentation:
\[
\mathcal{G}\cong\left\langle a,b,c\,:\,b^{-1}ab=a^{2},c^{-1}ac=a^{2},bc=cb\right\rangle
\]
If $K=\{1,2\}$ or $K=\{1,3\}$, let $\mathcal{N}:=\lllangle bc^{-1}\rrrangle\trianglelefteq\mathcal{G}$ and let $\pi:\mathcal{G}\rightarrow\mathcal{G}/\mathcal{N}$ be the canonical projection. It is easy to see that both $\pi\circ\nu_{\{1,2\}}$ and $\pi\circ\nu_{\{1,3\}}$ are isomorphisms and, hence, both $\nu_{\{1,2\}}$ and $\nu_{\{1,3\}}$ are injective. If $K=\{2,3\}$, let $\mathcal{N}:=\llangle a\rrangle\trianglelefteq\mathcal{G}$ instead and proceed analogously. Finally, if $|K|\leq 1$, then $K$ is contained in some $\widetilde{K}\subseteq\{1,2,3\}$ with $|\widetilde{K}|=2$. By construction of the colimit, we have $\nu_{K}=\nu_{\widetilde{K}}\circ\varphi_{K\widetilde{K}}$. Since both $\nu_{\widetilde{K}}$ and $\varphi_{K\widetilde{K}}$ are injective, their composition is injective, too. \hfill $\Box$%
\begin{proposition}
\label{prn:examplepart2}
The equation $\mathcal{G}|_{\{1,2\}}\cap\mathcal{G}|_{\{1,3\}}=\mathcal{G}|_{\{1\}}$ does not hold.
\end{proposition}
\textbf{Proof.} Use the isomorphism $\mu_{\{1,2\}}:G_{\{1,2\}}\rightarrow\mathcal{G}|_{\{1,2\}}$ and Britton's Lemma to show that the word $bab^{-1}$ represents an element in $\mathcal{G}|_{\{1,2\}}$ that is not in $\mathcal{G}|_{\{1\}}$. On the other hand, in the colimit group $\mathcal{G}$, the following equations hold: \[ bab^{-1}=bca^{2}c^{-1}b^{-1}=cba^{2}b^{-1}c^{-1}=cac^{-1} \] So, the words $bab^{-1}$ and $cac^{-1}$ represent the same element of the colimit group $\mathcal{G}$, which is in $\mathcal{G}|_{\{1,2\}}\cap\mathcal{G}|_{\{1,3\}}$ but not in $\mathcal{G}|_{\{1\}}$. The above calculation is also illustrated in Figure~\ref{fig:picture}, in terms of disc pictures. \hfill $\Box$
\subsection{Preliminaries about disc pictures} \label{sub:pictures}
\begin{figure}
\begin{center}
\input{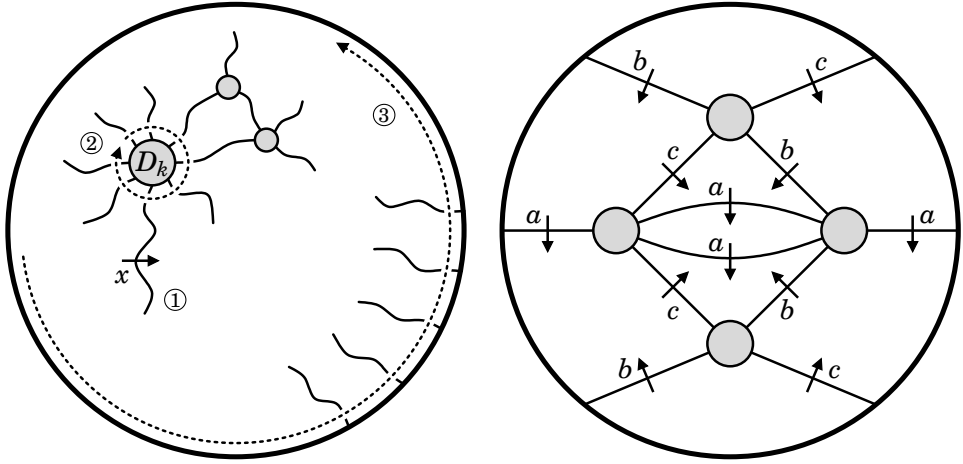}
\end{center}
\caption{Some elements of a disc picture (left) and an example showing that $bab^{-1}=cac^{-1}$ holds in $G\cong\left\langle a,b,c\,:\,b^{-1}ab=a^{2},c^{-1}ac=a^{2},bc=cb\right\rangle$ (right).}
\label{fig:picture}
\end{figure}
The proof of the intersection theorem involves \emph{disc pictures}. Let us therefore recall some basic notions from \cite{Corson1996}. Consider a group $G$ and a presentation $G\cong\langle X:R\rangle$. A disc picture $\mathcal{P}$ over this presentation consists of the disjoint union of closed discs $D_{1},D_{2},\dotsc,D_{n}$ in the interior of a closed disc $D$ and a compact 1-manifold $M$ properly embedded into $D\setminus\interior(D_{1}\cup D_{2}\cup\dotsb\cup D_{n})$.

\begin{figure}
\begin{center}
\input{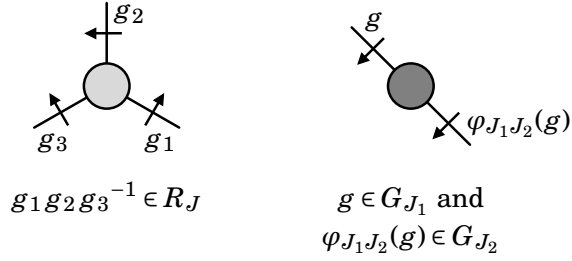}
\end{center}
\caption{A local vertex (bright) and a joining vertex (dark).}
\label{fig:vertices}
\end{figure}

The closed discs $D_{1},D_{2},\dotsc,D_{n}$ are called \emph{vertices}, the components of $M$ are called \emph{arcs}. Moreover, the components of $\interior(D)\setminus(D_{1}\cup D_{2}\cup\dotsb\cup D_{n}\cup M)$ are called \emph{regions}. Every arc has a transversal orientation and is labelled by a generator, see \circled{1} in Figure~\ref{fig:picture}. Every vertex $D_{k}$ has the property that one can read off a relator \emph{along} its boundary $\partial D_{k}$, i.\,e.~by starting at some point on $\partial D_{k}\setminus M$ and going once around $\partial D_{k}$ in some orientation, see \circled{2} in Figure~\ref{fig:picture}. Every word that can be read off along the outer boundary $\partial D$ is called a \emph{boundary word} of the disc picture, see \circled{3} in Figure~\ref{fig:picture}. It is well known, and easy to verify, that a word over the generators and their formal inverses represents the identity of the group $G$ if and only if it is a boundary word of some disc picture over the presentation $G\cong\langle X:R\rangle$. Disc pictures are, roughly speaking, duals of van Kampen diagrams. For further reading about them we refer to \cite{BogleyPride1993}. In addition to the above, the following notions will be of relevance for us.
\begin{definition}[``subpicture'']
\label{def:subpictures}
Consider a closed disc $D_{\mathcal{Q}}$ in $D$. If the parts of the disc picture $\mathcal{P}$ that are contained in $D_{\mathcal{Q}}$ assemble to a disc picture $\mathcal{Q}$, we call $\mathcal{Q}$ a subpicture of $\mathcal{P}$. Notice that every boundary word of a subpicture does necessarily represent the identity of the group $G$. A simple kind of subpicture is a spider. It consists of exactly one vertex $D_{k}$ and some arcs, each of which connects $D_{k}$ to the outer boundary $\partial D_{\mathcal{Q}}$ of the subpicture $\mathcal{Q}$.
\end{definition}
Since we are interested in Corson diagrams and their colimit groups, we will focus on disc pictures over (\ref{def:colimit}). Here, it makes sense to distinguish between local and joining vertices.
\begin{definition}[``local and joining vertices'']
\label{def:localandjoiningvertices}
A vertex $D_{k}$ is called local if one can read off a relator of the form ${g_{1}}^{\varepsilon_{1}}{g_{2}}^{\varepsilon_{2}}\dotsm{g_{m}}^{\varepsilon_{m}}\in R_{J}$ along its boundary. Otherwise, it is called joining, in which case one can read off a relator of the form $g=\varphi_{J_{1}J_{2}}(g)$ with $g\in G_{J_{1}}$ and $\varphi_{J_{1}J_{2}}(g)\in G_{J_{2}}$, see Figure~\ref{fig:vertices}.
\end{definition}
\begin{definition}[``bridge'']
\label{def:bridges}
Let $\mathcal{B}$ be the union of the compact 1-manifold $M$ and the joining vertices. The components of $\mathcal{B}$ are called bridges. Every simply connected bridge has two distinct endpoints, each of which lies either on the boundary of some local vertex or on the outer boundary. Two local vertices, say $D_{k}$ and $D_{l}$, are called neighbours if there is a bridge that connects $D_{k}$ and $D_{l}$, i.\,e.~a bridge with one endpoint on the boundary $\partial D_{k}$ and the other endpoint on the boundary $\partial D_{l}$.
\end{definition}
\begin{definition}[``inner and outer'']
\label{def:innerouter}
A bridge is called inner if it connects two local vertices. Similarly, a region is called inner if its closure does not meet the outer boundary $\partial D$. A bridge or a region that is not inner is called outer.
\end{definition}
Let us consider a non-spherical Corson diagram! As stated in Remark~\ref{rem:shorthand}, we may interpret $\mathcal{G}|_{K}$ as a subgroup of $\mathcal{G}$. The following lemma uses this interpretation to describe the labels of the arcs of a bridge.
\begin{lemma}
\label{lem:bridge}
Consider a bridge with $m$ arcs that are labelled by generators $b_{1}\in G_{J_{1}},b_{2}\in G_{J_{2}},\dotsc,b_{m}\in G_{J_{m}}$. Then, all these generators represent the same element of the colimit group $\mathcal{G}$. This element, say $b\in\mathcal{G}$, is called the value of the bridge. It is contained in the subgroup $\mathcal{G}|_{J_{1}\cap J_{2}\cap\dotsb\cap J_{m}}$.
\end{lemma}
\textbf{Proof.} The first assertion is immediate. So, we only need to verify that the value of the bridge is actually contained in $\mathcal{G}|_{J_{1}\cap J_{2}\cap\dotsb\cap J_{m}}$. Let us make two observations! First, if one of the sets $J_{1},J_{2},\dotsc,J_{m}$ is empty, say $J_{i}=\emptyset$, then the value of the bridge can be represented by $b_{i}\in G_{\emptyset}$. So, $b\in\mathcal{G}|_{\emptyset}$. This, of course, can be written as $b\in\mathcal{G}|_{J_{1}\cap J_{2}\cap\dotsb\cap J_{m}}$, whence we are done. Therefore, we may assume w.\,l.\,o.\,g.~that none of the sets $J_{1},J_{2},\dotsc,J_{m}$ is empty, in which case they must alternately have cardinality 1 and 2. Second, if $m=1$, there is nothing to show. So, we may assume w.\,l.\,o.\,g.~that $m\geq 2$. But then, there must be at least one set of cardinality 1 among $J_{1},J_{2},\dotsc,J_{m}$.
\begin{enumerate}
\item If all the sets of cardinality 1 are equal, say equal to $\{i\}$, then $b\in\mathcal{G}|_{\{i\}}$. But, in this situation, all sets of cardinality 2 must contain $\{i\}$ as a subset, which implies that $J_{1}\cap J_{2}\cap\dotsb\cap J_{m}=\{i\}$. Therefore, $b\in\mathcal{G}|_{\{i\}}$ can be written as $b\in\mathcal{G}|_{J_{1}\cap J_{2}\cap\dotsb\cap J_{m}}$.
\item If there are two distinct sets of cardinality 1 among $J_{1},J_{2},\dotsc,J_{m}$, say $\{i\}$ and $\{j\}$, then $b\in\mathcal{G}|_{\{i\}}\cap\mathcal{G}|_{\{j\}}$. We claim that $\mathcal{G}|_{\{i\}}\cap\mathcal{G}|_{\{j\}}=\mathcal{G}|_{\emptyset}$. Once this has been shown, we know that $b\in\mathcal{G}|_{\emptyset}$. Since $J_{1}\cap J_{2}\cap\dotsb\cap J_{m}=\emptyset$, this can, again, be written as $b\in\mathcal{G}|_{J_{1}\cap J_{2}\cap\dotsb\cap J_{m}}$.

It remains to show that $\mathcal{G}|_{\{i\}}\cap\mathcal{G}|_{\{j\}}=\mathcal{G}|_{\emptyset}$. By Remark~\ref{rem:anglesinterpretationpart1}, we already know that the equation $\varphi_{\{i\}\{i,j\}}\left(G_{\{i\}}\right)\cap\varphi_{\{j\}\{i,j\}}\left(G_{\{j\}}\right)=\varphi_{\emptyset\{i,j\}}\left(G_{\emptyset}\right)$ holds. In order to transport this equation to the colimit group $\mathcal{G}$, we apply the injective homomorphism $\nu_{\{i,j\}}=\tilde\nu_{\{i,j\}}\circ\mu_{\{i,j\}}:G_{\{i,j\}}\rightarrow\mathcal{G}$ and obtain:
\[
\nu_{\{i,j\}}\left(\varphi_{\{i\}\{i,j\}}\left(G_{\{i\}}\right)\right)\cap\nu_{\{i,j\}}\left(\varphi_{\{j\}\{i,j\}}\left(G_{\{j\}}\right)\right)=\nu_{\{i,j\}}\left(\varphi_{\emptyset\{i,j\}}\left(G_{\emptyset}\right)\right)
\]
Since the images $\nu_{\{i,j\}}\left(\varphi_{K\{i,j\}}\left(G_{K}\right)\right)$ with $K\subset\{i,j\}$ satisfy the equations $\nu_{\{i,j\}}\left(\varphi_{K\{i,j\}}\left(G_{K}\right)\right)=\nu_{K}\left(G_{K}\right)=\tilde\nu_{K}\left(\mu_{K}\left(G_{K}\right)\right)=\tilde\nu_{K}\left(\mathcal{G}|_{K}\right)$, we finally obtain $\tilde\nu_{\{i\}}\left(\mathcal{G}|_{\{i\}}\right)\cap\tilde\nu_{\{j\}}\left(\mathcal{G}|_{\{j\}}\right)=\tilde\nu_{\emptyset}\left(\mathcal{G}|_{\emptyset}\right)$, which reads as $\mathcal{G}|_{\{i\}}\cap\mathcal{G}|_{\{j\}}=\mathcal{G}|_{\emptyset}$ in the shorthand notation of Remark~\ref{rem:shorthand}.\hfill $\Box$
\end{enumerate}
\subsection{Statement and proof of the intersection theorem} \label{sub:intersection}
We are now ready to discuss the intersection theorem. The proof is based on ideas and techniques that go back to Gersten and Stallings in \cite{Stallings1991} and Corson in \cite{Corson1996}.
\begin{theorem}
\label{thm:intersection}
For every non-spherical Corson diagram and every two subsets $K_{1},K_{2}\subseteq I$ the equation $\mathcal{G}|_{K_{1}}\cap\mathcal{G}|_{K_{2}}=\mathcal{G}|_{K_{1}\cap K_{2}}$ holds. 
\end{theorem}
\textbf{Proof.} The inclusion ``$\supseteq$'' is a consequence of Remark~\ref{rem:shorthandinterpretation}. So, we only need to verify the inclusion ``$\subseteq$''. Suppose there were a non-spherical Corson diagram and two subsets $K_{1},K_{2}\subseteq I$ with $\mathcal{G}|_{K_{1}}\cap\mathcal{G}|_{K_{2}}\not\subseteq\mathcal{G}|_{K_{1}\cap K_{2}}$. Then, we can find an element $g\in\mathcal{G}$ with $g\in\mathcal{G}|_{K_{1}}\cap\mathcal{G}|_{K_{2}}$ but $g\not\in\mathcal{G}|_{K_{1}\cap K_{2}}$. This element can be represented by words $w_{1}$ and $w_{2}$ in the following languages:
\[
w_{1}\in\left(\bigsqcup_{\genfrac{}{}{0pt}{}{J\subseteq K_{1}}{|J|\leq 2}}G_{J}\sqcup{G_{J}}^{-}\right)^{\ast}\quad\text{and}\quad w_{2}\in\left(\bigsqcup_{\genfrac{}{}{0pt}{}{J\subseteq K_{2}}{|J|\leq 2}}G_{J}\sqcup{G_{J}}^{-}\right)^{\ast}
\]
Since $w_{1}$ and $w_{2}$ represent the same element of the colimit group $\mathcal{G}$, there is a disc picture $\mathcal{P}$ over (\ref{def:colimit}) with boundary word $w_{1}{w_{2}}^{-1}$. By construction, $g\not\in\mathfrak{G}|_{K_{1}\cap K_{2}}$. So, it cannot be the identity of the colimit group $\mathcal{G}$. Therefore, the words $w_{1}$ and $w_{2}$ cannot be empty and there are at least two arcs (or one arc twice) incident with the outer boundary $\partial D$.

We may assume w.\,l.\,o.\,g.~that the element $g$, the words $w_{1}$ and $w_{2}$, and the disc picture $\mathcal{P}$ are chosen in such a way that the \emph{complexity} of the disc picture is minimal, i.\,e.~the number of local vertices is minimal and, among all disc pictures with this minimal number of local vertices, the number of bridges is minimal. This assumption has many consequences on the structure of the disc picture.
\begin{enumerate}
\item \label{item:connected} \textbf{The disc picture $\mathcal{P}$ is connected. In particular, since there are arcs incident with the outer boundary $\partial D$, every local vertex is incident with at least one arc. Moreover, all bridges and regions are simply connected.} We claim that if $\mathcal{P}$ was not connected, we could remove at least one component and, hence, obtain a disc picture with fewer local vertices or with the same number of local vertices but fewer bridges. In other words, we could obtain a disc picture of lower complexity.

First, notice that there are two distinct points $x,y\in\partial D\setminus M$ such that one can read off the words $w_{1}$ and $w_{2}$ when going from $x$ to $y$ along the respective side of $\partial D$, see \circled{1} in Figure~\ref{fig:simplify-3}. If there is a component of $\mathcal{P}$ that is incident with at most one side of $\partial D$, we can remove it. In this case, the boundary words of the disc picture may change. But the new disc picture gives rise to new words $\widetilde{w}_{1}$ and $\widetilde{w}_{2}$. Since the removed component has been incident with at most one side of $\partial D$, at least one of the words $\widetilde{w}_{i}$ is equal to $w_{i}$. So, $\widetilde{w}_{1}$ and $\widetilde{w}_{2}$, which represent the same element of the colimit group $\mathcal{G}$, both still represent $g$.

\textit{In the following steps, as here, we may obtain new words $\widetilde{w}_{1}$ and $\widetilde{w}_{2}$, and sometimes even a new element $\widetilde{g}\in\mathcal{G}$. But, in each step, it is easy to see that this data could have been chosen right at the beginning.}

So, we may assume w.\,l.\,o.\,g.~that every component of $\mathcal{P}$ is incident with both sides of $\partial D$. Suppose there is more than one such component and let $C$ be the first one traversed when going from $x$ to $y$ along $\partial D$. For a moment, let us focus on $C$ and ignore all the other components! Now, one can read off new words $\widetilde{w}_{1}$ and $\widetilde{w}_{2}$ along the respective sides of $\partial D$ that represent a new element $\widetilde{g}\in\mathcal{G}$, see \circled{2} in Figure~\ref{fig:simplify-3}. By construction, $\widetilde{g}\in\mathcal{G}|_{K_{1}}\cap\mathcal{G}|_{K_{2}}$. If $\widetilde{g}\not\in\mathcal{G}|_{K_{1}\cap K_{2}}$, the component $C$ is already a suitable disc picture and we can remove the other components completely. On the other hand, if $\widetilde{g}\in\mathcal{G}|_{K_{1}\cap K_{2}}$, we can keep the other components and remove $C$. The words that one can now read off along the respective sides of $\partial D$ represent the element $\widetilde{\widetilde{g}}:=\widetilde{g}^{-1}g\in\mathcal{G}$, which is the product of an element in $\mathcal{G}|_{K_{1}\cap K_{2}}$ and an element not in $\mathcal{G}|_{K_{1}\cap K_{2}}$. Therefore, $\widetilde{\widetilde{g}}\not\in\mathcal{G}|_{K_{1}\cap K_{2}}$ and, again, we end up with a suitable disc picture.~\Checkmark
\begin{figure}
\begin{center}
\input{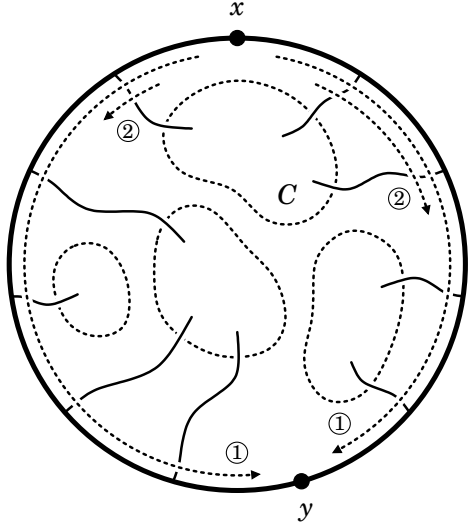}
\end{center}
\caption{The disc picture $\mathcal{P}$ is connected.}
\label{fig:simplify-3}
\end{figure}
\begin{definition}[``type of a local vertex'']
\label{def:type}
Since every local vertex is incident with at least one arc, we can associate a type to every local vertex. More precisely, for every local vertex $D_{k}$ there is a unique\footnote{Recall that the set of generators is the disjoint union of all $G_{J}$ with $J\subseteq I$ and $|J|\leq 2$.} subset $J\subseteq I$ with $|J|\leq 2$ such that all arcs incident with $D_{k}$ are labelled by generators from $G_{J}$. In this case, we say that $D_{k}$ is of type $J$. 
\end{definition}
\item \label{item:localvertex} \textbf{The disc picture $\mathcal{P}$ has at least one local vertex.} If $\mathcal{P}$ had no local vertex at all, it would have to be a single bridge $B$ connecting the two sides of $\partial D$.

Depending on the transversal orientation of its arcs, the value of $B$ is either $g$ or $g^{-1}$. The extremal arcs of $B$ are labelled by generators, say $b_{1}\in G_{J_{1}}$ and $b_{m}\in G_{J_{m}}$ with $J_{1}\subseteq K_{1}$ and $J_{m}\subseteq K_{2}$. Using Lemma~\ref{lem:bridge} and Remark~\ref{rem:shorthandinterpretation}, we can now observe that $g\in\mathcal{G}|_{J_{1}\cap J_{2}\cap\dotsb\cap J_{m}}\subseteq\mathcal{G}|_{J_{1}\cap J_{m}}\subseteq\mathcal{G}|_{K_{1}\cap K_{2}}$, in contradiction to $g\not\in\mathcal{G}|_{K_{1}\cap K_{2}}$.~\Checkmark
\item \label{item:noloop} \textbf{A local vertex cannot be neighbour of itself.} If there was such a local vertex $D_{k}$, say of type $J$, we could consider the subpicture $\mathcal{Q}$ consisting of the local vertex $D_{k}$, one of the bridges that connect $D_{k}$ with itself, everything that is enclosed by this bridge, and the extremal parts of the remaining arcs incident with $D_{k}$, see \circled{1} in Figure~\ref{fig:simplify-1-2}. Every boundary word $w$ of the subpicture $\mathcal{Q}$ is a word over generators from $G_{J}$ and their formal inverses that represents the identity of the colimit group $\mathcal{G}$. Since the natural homomorphism $\nu_{J}:G_{J}\rightarrow\mathcal{G}$ is injective, the word $w$ does not only represent the identity of the colimit group $\mathcal{G}$ but also the identity of the group $G_{J}$. Therefore, we can remove the subpicture $\mathcal{Q}$ and replace it by a single spider with boundary word $w$, see \circled{2} in Figure~\ref{fig:simplify-1-2}. After this modification, we obtain a disc picture with at most as many local vertices and stricly fewer bridges, and, hence, of lower complexity.~\Checkmark
\item \label{item:noneighbours} \textbf{Two local vertices of the same type cannot be neighbours.} If there were two such local vertices $D_{k}$ and $D_{l}$, w.\,l.\,o.\,g.~$D_{k}\neq D_{l}$, we could consider the subpicture $\mathcal{Q}$ consisting of the local vertices $D_{k}$ and $D_{l}$, one of the bridges that connect $D_{k}$ and $D_{l}$, and the extremal parts of the remaining arcs incident with $D_{k}$ and $D_{l}$, see \circled{3} in Figure~\ref{fig:simplify-1-2}. By the same arguments as in (\ref{item:noloop}), we can remove the subpicture $\mathcal{Q}$ and replace it by a single spider with the same boundary word, see \circled{4} in Figure~\ref{fig:simplify-1-2}. Again, we obtain a disc picture of lower complexity.~\Checkmark
\begin{figure}
\begin{center}
\input{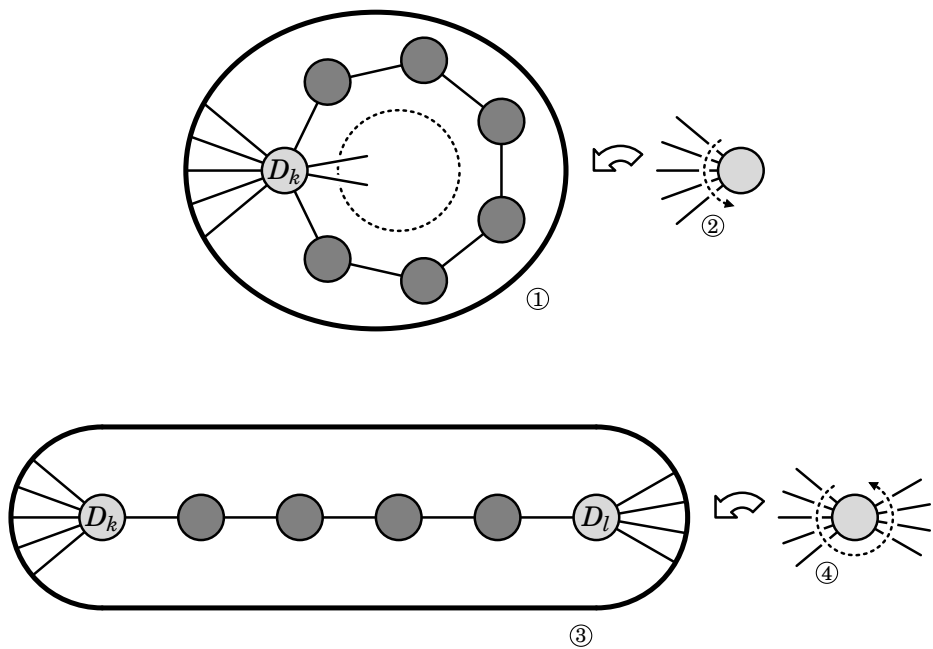}
\end{center}
\caption{Replace the subpictures $\mathcal{Q}$ by spiders.}
\label{fig:simplify-1-2}
\end{figure}
\item \label{item:noshortbridge} \textbf{Every bridge has at least two arcs.} If there was a bridge $B$ with only one arc, we would consider three cases: First, if $B$ is connecting two local vertices, say of types $J_{1}$ and $J_{2}$, then $J_{1}=J_{2}$, in contradiction to (\ref{item:noneighbours}). Second, if $B$ is connecting the outer boundary $\partial D$ with itself, then, by (\ref{item:connected}), $B$ is already the whole disc picture, in contradiction to (\ref{item:localvertex}). So, we may assume w.\,l.\,o.\,g.~that $B$ is connecting a local vertex $D_{k}$, say of type $J$, and the outer boundary $\partial D$, say at the side of $\partial D$ along which one can read off the word $w_{1}$.

By the former, the bridge $B$ is labelled by some generator from $G_{J}$ and, by the latter, $J\subseteq K_{1}$. Now, consider the subpicture $\mathcal{Q}$ consisting of the local vertex $D_{k}$, the bridge $B$, and the extremal parts of the remaining arcs incident with $D_{k}$, see \circled{1} in Figure~\ref{fig:simplify-8}. Replace it by a subpicture in which the arcs traversing $\partial D_{\mathcal{Q}}$, which are all labelled by generators from $G_{J}$, are extended to the outer boundary $\partial D$, see \circled{2} in Figure~\ref{fig:simplify-8}. This gives a disc picture with one fewer local vertex and, hence, of lower complexity.~\Checkmark
\begin{figure}
\begin{center}
\input{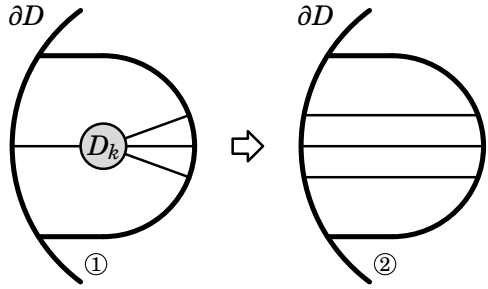}
\end{center}
\caption{Replace the subpicture $\mathcal{Q}$ by some arcs.}
\label{fig:simplify-8}
\end{figure}
\item \label{item:leftright} \textbf{The two regions on either side of a bridge cannot be the same.} Suppose there was a bridge $B$ having the same region $R$ on either side. Then, we can find a subpicture $\mathcal{Q}$ whose boundary $\partial D_{\mathcal{Q}}$ is contained in $R$, except of one point where it crosses $B$, see Figure~\ref{fig:simplify-4}. Therefore, the value of $B$ is the identity of the colimit group $\mathcal{G}$. Since the natural homomorphisms $\nu_{J}:G_{J}\rightarrow\mathcal{G}$ are injective, the labels of the arcs of $B$ must also be the identities of the respective groups $G_{J}$. So, we can remove the bridge $B$ and obtain a disc picture of lower complexity.~\Checkmark
\begin{figure}
\begin{center}
\input{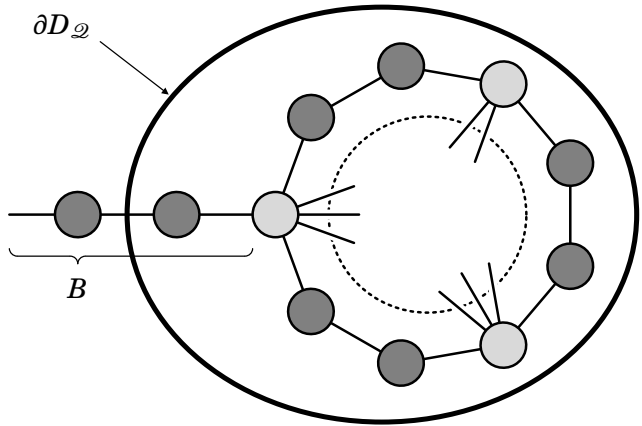}
\end{center}
\caption{The value of $B$ is the identity.}
\label{fig:simplify-4}
\end{figure}
\item \label{item:value} \textbf{The value of a bridge cannot be an element of $\mathcal{G}|_{\emptyset}$.} Suppose there was a bridge $B$ with value $b\in\mathcal{G}|_{\emptyset}$. Since $b\in\mathcal{G}|_{\emptyset}$, there is a generator $b_{\emptyset}\in G_{\emptyset}$ that represents it. In fact, for every $J\subseteq I$ with $1\leq|J|\leq 2$ there is a generator $b_{J}\in G_{J}$ that represents it, namely $b_{J}:=\varphi_{\emptyset J}(b_{\emptyset})\in G_{J}$.

By (\ref{item:leftright}), the two regions on either side of $B$ cannot be the same. Among these two distinct regions choose the region $R$ with the property that the arcs of $B$ are heading away from $R$. Now, remove $B$ and relabel all the remaining arcs in the boundary $\partial R$ as follows: If an arc is labelled by a generator $a\in G_{J}$ and is heading towards $R$, relabel it by $ab_{J}\in G_{J}$. If it is heading away from $R$, relabel it by ${b_{J}}^{-1}a\in G_{J}$. This guarantees that one can still read off relators along the boundaries of the remaining vertices. Here, we leave the details to the reader, see \cite[Appendix]{Corson1996} for another description and \circled{1}--\circled{3} in Figure~\ref{fig:simplify-5} for some examples.

Notice that, in \circled{2} in Figure~\ref{fig:simplify-5}, the generator $\varphi_{J_{1}J_{2}}(a)\cdot b_{J_{2}}\in G_{J_{2}}$ satisfies the following equation:
\[ \begin{array}{r@{\;}c@{\;}l} \varphi_{J_{1}J_{2}}(a)\cdot b_{J_{2}} & = & \varphi_{J_{1}J_{2}}(a)\cdot\varphi_{\emptyset J_{2}}(b_{\emptyset}) \smallskip \\ & = & \left\{ \begin{array}{ll} \varphi_{J_{1}J_{2}}(a)\cdot\varphi_{J_{1}J_{2}}(b_{J_{1}}) & \text{if}~J_{1}=\emptyset \\ \varphi_{J_{1}J_{2}}(a)\cdot\varphi_{J_{1}J_{2}}(\varphi_{\emptyset J_{1}}(b_{\emptyset})) & \text{otherwise} \end{array} \right. \smallskip \\ & = & \varphi_{J_{1}J_{2}}(a)\cdot\varphi_{J_{1}J_{2}}(b_{J_{1}}) \smallskip \\ & = & \varphi_{J_{1}J_{2}}(a\cdot b_{J_{1}}) \end{array} \]
Therefore, one can actually read off the relator $a\cdot b_{J_{1}}=\varphi_{J_{1}J_{2}}(a\cdot b_{J_{1}})$ along the boundary of the respective joining vertex.
\begin{figure}
\begin{center}
\input{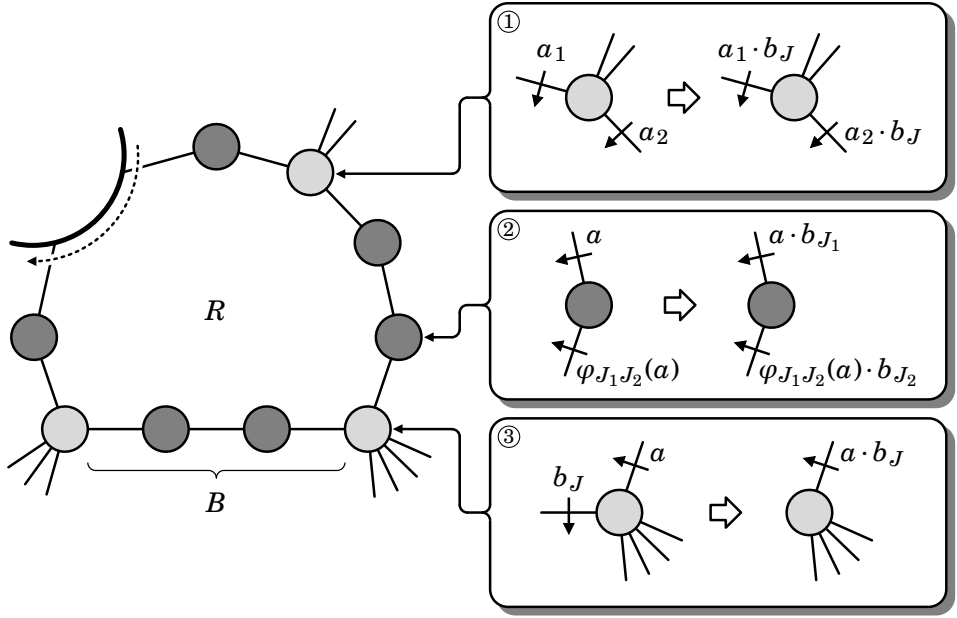}
\end{center}
\caption{Relabel the remaining arcs in the boundary $\partial R$.}
\label{fig:simplify-5}
\end{figure}

So, we obtain a disc picture with the same number of local vertices and strictly fewer bridges and, hence, of lower complexity.~\Checkmark
\item \label{item:neueins} \textbf{For every bridge there is a unique element $i\in I$ such that the value of the bridge is in $\mathcal{G}|_{\{i\}}\setminus\mathcal{G}|_{\emptyset}$.} Let $B$ be a bridge with $m$ arcs that are labelled by generators $b_{1}\in G_{J_{1}},b_{2}\in G_{J_{2}},\dotsc,b_{m}\in G_{J_{m}}$. By (\ref{item:value}), none of the sets $J_{1},J_{2},\dotsc,J_{m}$ is empty, which implies that they must alternately have cardinality 1 and 2. By (\ref{item:noshortbridge}), every bridge has at least two arcs, so $m\geq 2$. Therefore, there must be at least one set $J_{k}$ of cardinality 1, say $J_{k}=\{i\}$. So, the value of $B$ is an element of $\mathcal{G}|_{\{i\}}$ and, again by (\ref{item:value}), cannot be an element of $\mathcal{G}|_{\emptyset}$. The element $i\in I$ is unique: In the second step of the proof of Lemma~\ref{lem:bridge}, we have seen that $\mathcal{G}|_{\{i\}}\cap\mathcal{G}|_{\{j\}}=\mathcal{G}|_{\emptyset}$ for any two distinct $i,j\in I$. So, $B$ cannot have a value that is simultaneously contained in $\mathcal{G}|_{\{i\}}\setminus\mathcal{G}|_{\emptyset}$ and in $\mathcal{G}|_{\{j\}}\setminus\mathcal{G}|_{\emptyset}$.~\Checkmark
\begin{definition}[``type of a bridge'']
In this situation, we say that the bridge is of type $i$.
\end{definition}
\item \label{item:neuzwei} \textbf{If a bridge of type $i$ is incident with a local vertex of type $J$, then $i\in J$. Similarly, if it is incident with one side of the outer boundary $\partial D$, then $i\in K_{1}$ or $i\in K_{2}$ respectively.} We give a proof of the first assertion, the proof of the second one is essentially the same. Let $B$ be a bridge of type $i$ and let $D_{k}$ be a local vertex of type $J$. As we have seen in (\ref{item:neueins}), one of the arcs of $B$ is labelled by a generator from $G_{\{i\}}$. Moreover, the extremal arc of $B$ that is incident with $D_{k}$ is labelled by a generator from $G_{J}$. By Lemma~\ref{lem:bridge} and Remark~\ref{rem:shorthandinterpretation}, we can now observe that the bridge has a value in $\mathcal{G}|_{\{i\}\cap J}$. But, by (\ref{item:value}), the value is not in $\mathcal{G}|_{\emptyset}$. Therefore, $\{i\}\cap J\neq\emptyset$, whence $i\in J$.~\Checkmark
\item \label{item:nozerotype} \textbf{There are no local vertices of type $\emptyset$.} By (\ref{item:connected}), every local vertex is incident with at least one arc. So, if there was a local vertex of type $\emptyset$, it would have to be incident with an arc that is labelled by some generator $a\in G_{\emptyset}$. But this arc is part of a bridge with value in $\mathcal{G}|_{\emptyset}$, in contradiction to (\ref{item:value}).~\Checkmark
\item \label{item:noonetype} \textbf{There are no local vertices of type $\{i\}$ with $i\in I$.} First, observe that if there was such a local vertex $D_{k}$, it would have to be neighbour of some other local vertex $D_{l}$.

Suppose it wasn't. Then, all bridges that are incident with $D_{k}$ must either connect it to itself, which is not possible by (\ref{item:noloop}), or to the outer boundary $\partial D$. But, by (\ref{item:connected}), the disc picture $\mathcal{P}$ is connected. So, this is already the whole disc picture. In particular, all bridges are incident with $D_{k}$, which is a local vertex of type $\{i\}$. Therefore, by (\ref{item:neuzwei}), all bridges are of type $i$. But this means that each letter of $w_{1}$ and $w_{2}$ represents an element in $\mathcal{G}|_{\{i\}}$, whence $g\in\mathcal{G}|_{\{i\}}$.

On the other hand, since both $w_{1}$ and $w_{2}$ are not empty, there is at least one bridge connecting $D_{k}$ to either side of $\partial D$. Again, by (\ref{item:neuzwei}), this implies both $i\in K_{1}$ and $i\in K_{2}$. But since $\{i\}\subseteq K_{1}\cap K_{2}$, we can use Remark~\ref{rem:shorthandinterpretation} to conclude that $g\in\mathcal{G}|_{K_{1}\cap K_{2}}$, in contradiction to $g\not\in\mathcal{G}|_{K_{1}\cap K_{2}}$.

So, we may assume w.\,l.\,o.\,g.~that $D_{k}$ is neighbour of some other local vertex $D_{l}$. Consider a bridge that connects $D_{k}$ and $D_{l}$. Again, by (\ref{item:neuzwei}), this bridge is of type $i$ and the local vertex $D_{l}$ is of some type $J$ with $i\in J$, i.\,e.~$\{i\}\subseteq J$. By (\ref{item:noneighbours}), two local vertices of the same type cannot be neighbours. So, we actually obtain that $\{i\}\subset J$. Now, replace every arc that is incident with $D_{k}$, say labelled by some generator $a\in G_{\{i\}}$, by a sequence of three arcs with the same transversal orientation. The first and the third are labelled by $a\in G_{\{i\}}$, the second by $\varphi_{\{i\}J}(a)\in G_{J}$, see \circled{1} in Figure~\ref{fig:simplify-7}. Then, consider the subpicture $\mathcal{Q}$ indicated in \circled{2} in Figure~\ref{fig:simplify-7}. By the same arguments as in (\ref{item:noloop}) and (\ref{item:noneighbours}), we can remove the subpicture $\mathcal{Q}$ and replace it by a single spider with the same boundary word. Again, we obtain a disc picture of lower complexity.~\Checkmark
\begin{figure}
\begin{center}
\input{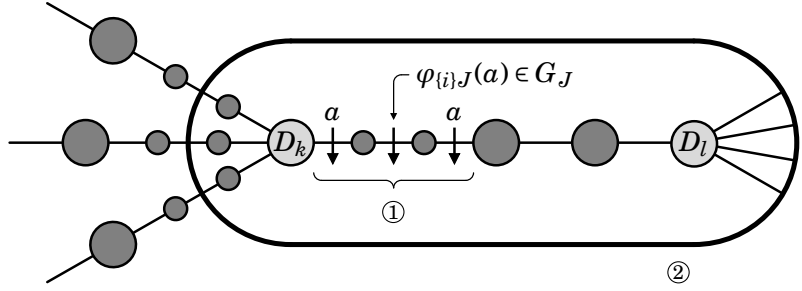}
\end{center}
\caption{Replace each arc by a sequence of three arcs.}
\label{fig:simplify-7}
\end{figure}
\begin{definition}[``angle and swap'']
By (\ref{item:nozerotype}) and (\ref{item:noonetype}), we know that every local vertex $D_{k}$ is of some type $\{i,j\}$ with distinct $i,j\in I$. By (\ref{item:neuzwei}), such a local vertex is incident with bridges each of which is either of type $i$ or of type $j$. We will consider the connected components of $\partial D_{k}\setminus M$. They are called angles. Moreover, an angle is called a swap if one of the two bridges enclosing it is of type $i$ and the other one is of type $j$.
\end{definition}
\item \label{item:nozeroswap} \textbf{Every local vertex has at least one swap in its boundary.} Suppose there was a local vertex $D_{k}$, say as above of type $\{i,j\}$, without any swap in its boundary. Then, all bridges that are incident with $D_{k}$ are of the same type, say of type $i$.

Since $D_{k}$ is a local vertex of type $\{i,j\}$, the arcs that are incident with $D_{k}$ are labelled by generators $a_{1},a_{2},\dotsc,a_{m}\in G_{\{i,j\}}$. The respective bridges are all of type $i$. So, each of these generators represents an element in $\mathfrak{G}|_{\{i\}}\smallsetminus\mathfrak{G}|_{\varnothing}$, whence we can even find generators $\widetilde{a}_{1},\widetilde{a}_{2},\dotsc,\widetilde{a}_{m}\in G_{\{i\}}$ representing the same elements, i.\,e.~satisfying the equations: \[ \nu_{\{i,j\}}\big(a_{s}\big)=\nu_{\{i\}}\big(\widetilde{a}_{s}\big)=\nu_{\{i,j\}}\left(\varphi_{\{i\}\{i,j\}}\left(\widetilde{a}_{s}\right)\right) \] Now, we can use the injectivity of the homomorphism $\nu_{\{i,j\}}:G_{\{i,j\}}\rightarrow\mathfrak{G}$ to conlcude that $a_{s}=\varphi_{\{i\}\{i,j\}}\left(\widetilde{a}_{s}\right)$.

Similarly to the modification described above in (\ref{item:noonetype}), we replace every arc that is incident with $D_{k}$ by a sequence of two arcs with the same transversal orientation. If the arc has been labelled by $a_{s}\in G_{\{i,j\}}$, the new arc that is incident with $D_{k}$ is labelled by $\widetilde{a}_{s}\in G_{\{i\}}$ whereas the other one is labelled by $a_{s}\in G_{\{i,j\}}$, see Figure~\ref{fig:simplify-9}.

Now, let $\varepsilon_{1},\varepsilon_{2},\dotsc,\varepsilon_{m}\in\{-1,1\}$ such that the word ${a_{1}}^{\varepsilon_{1}}{a_{2}}^{\varepsilon_{2}}\dotsb{a_{m}}^{\varepsilon_{m}}$ has originally been a boundary word of $D_{k}$. After this modification, one can read off the word  $\text{$\widetilde{a}_{1}$}^{\varepsilon_{1}}\text{$\widetilde{a}_{2}$}^{\varepsilon_{2}}\dotsb\text{$\widetilde{a}_{m}$}^{\varepsilon_{m}}$ along $\partial D_{k}$. But: %
\[
\varphi_{\{i\}\{i,j\}}\left(\text{$\widetilde{a}_{1}$}^{\varepsilon_{1}}\text{$\widetilde{a}_{2}$}^{\varepsilon_{2}}\dotsb\text{$\widetilde{a}_{m}$}^{\varepsilon_{m}}\right)={a_{1}}^{\varepsilon_{1}}{a_{2}}^{\varepsilon_{2}}\dotsb{a_{m}}^{\varepsilon_{m}}=1~\text{in}~G_{\{i,j\}}
\]
The injectivity of $\varphi_{\{i\}\{i,j\}}:G_{\{i\}}\rightarrow G_{\{i,j\}}$ implies that $D_{k}$ has become a local vertex of type $\{i\}$ and can be removed as in (\ref{item:noonetype}). So, we obtain a disc picture of lower complexity.~\Checkmark
\begin{figure}
\begin{center}
\input{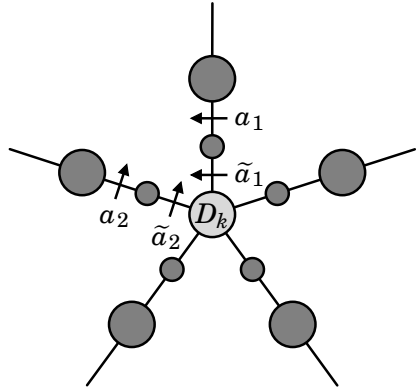}
\end{center}
\caption{Replace each arc by a sequence of two arcs.}
\label{fig:simplify-9}
\end{figure}
\setcounter{recallenumi}{\value{enumi}}
\end{enumerate}
We will now give a lower bound on the number of swaps in the boundary of a local vertex. So, assume given a local vertex $D_{k}$ of some type $\{i,j\}$ with distinct $i,j\in I$. We say that $D_{k}$ has \emph{sufficiently many swaps in its boundary} if the Gersten-Stallings angle $\sphericalangle_{\{i,j\}}\neq 0$ and the number of swaps $m\geq\nicefrac{2\pi\,}{\sphericalangle_{\{i,j\}}}$.
\begin{enumerate}
\setcounter{enumi}{\value{recallenumi}}
\item \label{item:swap} \textbf{Every local vertex has sufficiently many swaps in its boundary.} Suppose there was a local vertex $D_{k}$ without sufficiently many swaps in its boundary. By (\ref{item:nozeroswap}), there is at least one swap in its boundary. If we start at some swap and go from swap to swap once around $\partial D_{k}$, we can read off words $v_{1},v_{2},\dotsc,v_{m}$, see \circled{1} in Figure~\ref{fig:simplify-10}.

Each of these words represents an element in $\mathcal{G}$, so it makes sense to write $v_{1},v_{2},\dotsc,v_{m}\in\mathcal{G}$. Their product $v_{1}\cdot v_{2}\cdot\dotsc \cdot v_{m}=1$. We may assume w.\,l.\,o.\,g.~that $v_{1},v_{3},\dotsc,v_{m-1}\in\mathcal{G}|_{\{i\}}$ and $v_{2},v_{4},\dotsc,v_{m}\in\mathcal{G}|_{\{j\}}$. In order to show that at least one of these elements is contained in $\mathcal{G}|_{\emptyset}$, we use the injective homomorphisms $\nu_{K}:G_{K}\rightarrow\mathcal{G}$. The preimages:
\[ \begin{array}{l} {\nu_{\{i\}}}^{-1}(v_{1}),{\nu_{\{i\}}}^{-1}(v_{3}),\dotsc,{\nu_{\{i\}}}^{-1}(v_{m-1})\in G_{\{i\}} \smallskip \\ {\nu_{\{j\}}}^{-1}(v_{2}),{\nu_{\{j\}}}^{-1}(v_{4}),\dotsc,{\nu_{\{j\}}}^{-1}(v_{m})\in G_{\{j\}} \end{array} \]
assemble to an element $x:={\nu_{\{i\}}}^{-1}(v_{1})\cdot{\nu_{\{j\}}}^{-1}(v_{2})\cdot\dotsc\cdot{\nu_{\{j\}}}^{-1}(v_{m})$ of the amalgamated free product $G_{\{i\}}\ast_{G_{\varnothing}}G_{\{j\}}$. Now, recall the homomorphism $\alpha:G_{\{i\}}\ast_{G_{\varnothing}}G_{\{j\}}\rightarrow G_{\{i,j\}}$ introduced in Section~\ref{sub:curvature} and observe that:
\[ \begin{array}{r@{\;}c@{\;}l} \alpha(x) & = & \varphi_{\{i\}\{i,j\}}\left({\nu_{\{i\}}}^{-1}(v_{1})\right)\cdot\varphi_{\{j\}\{i,j\}}\left({\nu_{\{j\}}}^{-1}(v_{2})\right)\cdot\dotsc\cdot\varphi_{\{j\}\{i,j\}}\left({\nu_{\{j\}}}^{-1}(v_{m})\right) \smallskip \\ & = & {\nu_{\{i,j\}}}^{-1}(v_{1})\cdot{\nu_{\{i,j\}}}^{-1}(v_{2})\cdot\dotsc\cdot{\nu_{\{i,j\}}}^{-1}(v_{m}) \smallskip \\ & = & {\nu_{\{i,j\}}}^{-1}(v_{1}\cdot v_{2}\cdot\dotsc\cdot v_{m}) \smallskip \\ & = & {\nu_{\{i,j\}}}^{-1}(1) \smallskip \\ & = & 1 \end{array} \]
So, $x\in\operatorname{ker}(\alpha)$. Since $D_{k}$ does not have sufficiently many swaps in its boundary, we know that either the Gersten-Stallings angle $\sphericalangle_{\{i,j\}}=0$ or the length of $x$, which is at most $m$, is strictly smaller than $\nicefrac{2\pi\,}{\sphericalangle_{\{i,j\}}}$, which is nothing but the minimal length of a non-trivial element in $\operatorname{ker}(\alpha)$. In either case, $x$ must be trivial in $G_{\{i\}}\ast_{G_{\varnothing}}G_{\{j\}}$.

It is a consequence of the normal form theorem, see \cite[Lemma 1]{Miller1968}, that there is an index $k\in\{1,2,\dotsc,m\}$ such that ${\nu_{\{i\}}}^{-1}(v_{k})\in\varphi_{\emptyset\{i\}}(G_{\emptyset})$ or ${\nu_{\{j\}}}^{-1}(v_{k})\in\varphi_{\emptyset\{j\}}(G_{\emptyset})$, depending on the parity of $k$. But then:
\[ v_{k}\in\left\{\begin{array}{ll} \nu_{\{i\}}\left(\varphi_{\emptyset\{i\}}(G_{\emptyset})\right) & \text{if $k$ is odd} \\ \nu_{\{j\}}\left(\varphi_{\emptyset\{j\}}(G_{\emptyset})\right) & \text{if $k$ is even} \end{array}\right\}=\nu_{\emptyset}(G_{\emptyset})=\mathcal{G}|_{\emptyset} \]
In either case, $v_{k}\in\mathcal{G}|_{\emptyset}$. Now, we add a new local vertex and a new bridge to the disc picture as illustrated in \circled{2} in Figure~\ref{fig:simplify-10}: The arcs that had been traversed when reading off the word $v_{k}$ end up at the new local vertex, which is connected to $D_{k}$ by a single arc labelled by ${\nu_{\{i,j\}}}^{-1}(v_{k})$.

This increases both the number of local vertices and the number of bridges by 1. But still, some of the properties we have discussed so far, in particular (\ref{item:connected}), hold true and the bridge connecting the new local vertex and $D_{k}$, which has a value in $\mathcal{G}|_{\emptyset}$, can be removed as in (\ref{item:leftright}) or (\ref{item:value}). 

Next, we want to get rid of the new local vertex. If the removement of the bridge has made the disc picture $\mathcal{P}$ disconnected, we can remove one of the two components as in (\ref{item:connected}). Otherwise, there is still a path from the new local vertex to $D_{k}$, which implies that the new local vertex is neighbour of some other local vertex. Once this is clear, we can remove it as in (\ref{item:nozeroswap}) and in the final step of (\ref{item:noonetype}). In either case, in particular in the latter, we do not only remove the new local vertex but also at least one more bridge. So, in either case, we obtain a disc picture of lower complexity.~\Checkmark
\begin{figure}
\begin{center}
\input{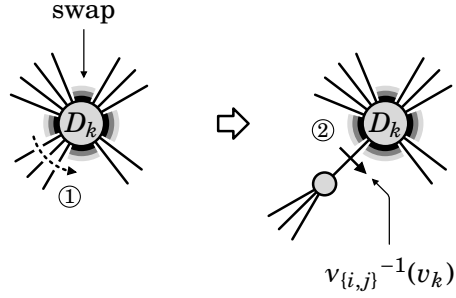}
\end{center}
\caption{Create a new local vertex and a new bridge.}
\label{fig:simplify-10}
\end{figure}
%
% Counter einprägen:
% \setcounter{recallenumi}{\value{enumi}}
%
\end{enumerate}
After all these observations, we can give an easy combinatorial argument that yields a contradiction. The principal idea is to dristibute weights over certain parts of the disc picture: Every local vertex gets the weight $2\pi$, every inner bridge gets the weight $-2\pi$, and every inner region gets the weight $2\pi$. For the notion of \emph{inner} see Definition~\ref{def:innerouter}. The weighted parts of the disc picture correspond to vertices, edges, and bounded regions of a planar graph, which is non-empty, finite, and connected. So, we may use Euler's formula for planar graphs to calculate the total weight:
\[ 2\pi\cdot\,\#\,\text{local vertices}-2\pi\cdot\,\#\,\text{inner bridges}+2\pi\cdot\,\#\,\text{inner regions}=2\pi \] 
Let us count again! But, this time, we reallocate the weights to the regions. Every inner bridge distributes its weight $-2\pi$ equally to the two regions on either side and every local vertex distributes its weight $2\pi$ equally to the swaps in its boundary, each of which lets it traverse to the adjacent region. The new total weight of a region $R$ is denoted by $\wt(R)$ and can be estimated from above using (\ref{item:traversingweight}), (\ref{item:nobadinnerregions}), and (\ref{item:fewbadregions}):
\begin{enumerate}[label=\textbf{(\alph*)},ref=\alph*]
%
% Counter auskramen:
% \setcounter{enumi}{\value{recallenumi}}
%
\item \label{item:traversingweight} \textbf{The weight traversing each swap is bounded by the Gersten-Stallings angle of the local vertex. In particular, it is at most $\nicefrac{\pi\,}{2}$.} Let $D_{k}$ be a local vertex of type $J=\{i,j\}$. By (\ref{item:swap}), $D_{k}$ has sufficiently many swaps in its boundary. So, the Gersten-Stallings angle $\sphericalangle_{\{i,j\}}\neq 0$ and the number of swaps is at least $\nicefrac{2\pi\,}{\sphericalangle_{\{i,j\}}}$. Since $D_{k}$ distributes its weight $2\pi$ equally to the swaps in its boundary, the weight traversing each swap is at most $\sphericalangle_{\{i,j\}}$.~\Checkmark
\item \label{item:nobadinnerregions} \textbf{There are no inner regions of positive weight.} Let $R$ be an inner region. By (\ref{item:connected}), $R$ is an open disc. The boundary $\partial R$ contains some number of inner bridges, say $m$, and the same number of angles, some of which may be swaps. By (\ref{item:noloop}), $m\geq 2$.

By (\ref{item:leftright}), each of the $m$ inner bridges contributes $-\pi$ to $\wt(R)$ and, by (\ref{item:traversingweight}), each of the at most $m$ swaps contributes at most $\nicefrac{\pi\,}{2}$. Therefore, we can estimate $\wt(R)$ as follows.

If $m\geq 4$, then $\wt(R)\leq 1\cdot2\pi-m\cdot\pi+m\cdot\nicefrac{\pi\,}{2}\leq 0$. If $m=3$ and there are at most two swaps, then $\wt(R)\leq 1\cdot 2\pi-3\cdot\pi+2\cdot\nicefrac{\pi\,}{2}=0$. If $m=3$ and there are exactly three swaps, then there are three pairwise distinct $i,j,k\in I$ such that the local vertices in the boundary $\partial R$ are of types $\{i,j\}$, $\{i,k\}$, $\{j,k\}$. Since we consider a non-spherical Corson diagram, there are no spherical triples. In particular, $\sphericalangle_{\{i,j\}}+\sphericalangle_{\{i,k\}}+\sphericalangle_{\{j,k\}}\leq\pi$, whence $\wt(R)\leq 1\cdot 2\pi-3\cdot\pi+\sphericalangle_{\{i,j\}}+\sphericalangle_{\{i,k\}}+\sphericalangle_{\{j,k\}}\leq 0$.

What remains is the case that $m=2$. If there were two distinct $i,j\in I$ such that one of the inner bridges is of type $i$ and the other is of type $j$, then both local vertices must be of type $\{i,j\}$, in contradiction to (\ref{item:noneighbours}). So, there cannot be any swap, whence $\wt(R)=1\cdot 2\pi-2\cdot\pi+0=0$.~\Checkmark
\item \label{item:fewbadregions} \textbf{There are at most two outer regions of positive weight. Each of them has at most weight $\nicefrac{\pi\,}{2}$.} Let $R$ be an outer region. By (\ref{item:connected}), $R$ is an open disc. The boundary $\partial R$ contains some number of bridges, say $m$, and some number of angles. By (\ref{item:connected}), (\ref{item:localvertex}), and (\ref{item:leftright}), exactly two of the bridges are outer, which implies that $m\geq 2$. Moreover, it contains exactly $m-1$ angles, some of which may be swaps.

If $m\geq 3$, then $\wt(R)\leq 0\cdot2\pi-(m-2)\cdot\pi+(m-1)\cdot\nicefrac{\pi\,}{2}=(-m+3)\cdot\nicefrac{\pi\,}{2}\leq 0$. If $m=2$ and there is no swap, then $\wt(R)=0\cdot 2\pi-0\cdot\pi+0=0$. What remains is the case that $m=2$ and there is a swap. Then, $\wt(R)$ might be positive, but $\wt(R)\leq 0\cdot 2\pi-0\cdot\pi+1\cdot\nicefrac{\pi\,}{2}=\nicefrac{\pi\,}{2}$.

Next, we show that this can happen at most twice. Since there is a swap, we know that there are two distinct $i,j\in I$ such that one of the outer bridges is of type $i$ and the other one is of type $j$. If both bridges end up at the same side of $\partial D$, say at the side of $\partial D$ along which one can read off the word $w_{1}$, then, by (\ref{item:neuzwei}), both $i\in K_{1}$ and $j\in K_{1}$. This allows us to remove the respective local vertex $D_{k}$ of type $\{i,j\}\subseteq K_{1}$ as in the final step of (\ref{item:noshortbridge}) and to obtain a disc picture of lower complexity. So, we may assume w.\,l.\,o.\,g.~that the two bridges end up at different sides of $\partial D$. But, by (\ref{item:connected}), we know that the disc picture $\mathcal{P}$ is connected. Hence, this can happen at most twice, namely when the boundary $\partial R$ contains one of the two points $x$ and $y$ that have been chosen in (\ref{item:connected}).~\Checkmark
\end{enumerate}
By (\ref{item:nobadinnerregions}) and (\ref{item:fewbadregions}), the total weight given to the disc picture is at most $\pi$. This is a contradiction to the above observation that the total weight amounts to $2\pi$, which completes the proof. \hfill $\Box$
\subsection{Interpretation} \label{sub:interpretation}
\begin{figure}
\begin{center}
\input{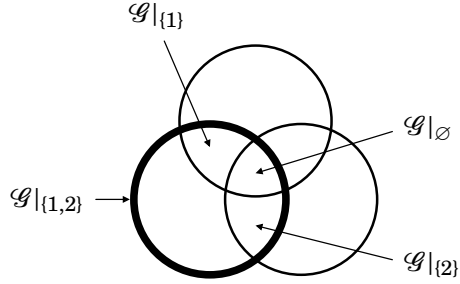}
\end{center}
\caption{$\mathcal{G}|_{\{1,2\}}$ intersecting $\mathcal{G}|_{\{1,3\}}$ and $\mathcal{G}|_{\{2,3\}}$.}
\label{fig:magma}
\end{figure}

In case of a non-spherical triangle of groups, the intersection theorem says that the groups $\mathcal{G}|_{K}$ with $K\subseteq\{1,2,3\}$ and $|K|\leq 2$ intersect exactly as sketched in Figure~\ref{fig:magma}. One particularly nice way of reading the intersection theorem is to start with such a setting: Let $M$ be the union of three groups with the property that each two of them intersect along a common subgroup, by which we implicitly mean that the two multiplications agree on the subgroup. $M$ is a set equipped with a partial multiplication and the question arises whether $M$ can be homomorphically embedded into a group, i.\,e.~whether there exists an injective map into a group such that the restriction to each of the three groups is a homomorphism. In order to give a partial answer to this question, we may interpret our three groups and their intersections, equipped with the inclusion maps, as a triangle of groups. By construction of the colimit, the natural homomorphisms $\nu_{J}:G_{J}\rightarrow\mathcal{G}$ agree on the intersections and, hence, yield a map $\nu:M\rightarrow\mathcal{G}$. It turns out that this map is injective if and only if the natural homomorphisms $\nu_{J}:G_{J}\rightarrow\mathcal{G}$ are injective (embedding theorem) and the equations $\mathcal{G}|_{K_{1}}\cap\mathcal{G}|_{K_{2}}=\mathcal{G}|_{K_{1}\cap K_{2}}$ hold (intersection theorem). So, if the triangle of groups is non-spherical, the answer is affirmative.
\begin{remark}
On the other hand, it is a consequence of the universal property, see \cite[Chapter 3]{AdamekHerrlichStrecker1990}, that if the map $\nu:M\rightarrow\mathcal{G}$ is not injective, then $M$ cannot be homomorphically embedded into any group and the answer is negative.
\end{remark}
%
% ==================================================
% §4 Billiards theorem for triangles of groups
% ==================================================
%
\section{Billiards theorem for triangles of groups} \label{sec:billiards}
In the previous section, we used a combinatorial argument based on Euler's formula for planar graphs to prove the intersection theorem. This kind of argument, be it in the language of homotopies, see e.\,g.~\cite{Stallings1991}, in the language of van Kampen diagrams, see e.\,g.~\cite{EdjvetRosenbergerStilleThomas2000} and \cite{HowieKopteva2006}, or in the language of disc pictures, see e.\,g.~\cite{Corson1996}, turned out to be very powerful in our context, and we should not miss to highlight the following two results from \cite{EdjvetRosenbergerStilleThomas2000} and \cite{HowieKopteva2006}. In each of them, one considers a non-spherical triangle of groups and assumes that for every $a\in\{1,2,3\}$ there is an element $g_{a}\in\mathcal{G}|_{\{a\}}\setminus\mathcal{G}|_{\emptyset}$.
\begin{theorem}[Edjvet et al.]
\label{thm:erst}
The product $g_{1}g_{2}g_{3}\in\mathcal{G}$ has infinite order. \hfill $\Box$
\end{theorem}
\begin{theorem}[Howie-Kopteva]
\label{thm:hk}
If the triangle of groups is hyperbolic, then there is an $n\in\mathbb{N}$ such that the elements $(g_{1}g_{2}g_{3})^{n}\in\mathcal{G}$ and $(g_{1}g_{3}g_{2})^{n}\in\mathcal{G}$ generate a non-abelian free subgroup. \hfill $\Box$
\end{theorem}
\begin{remark}
\label{rem:triviallyfilled}
Both, in Theorems~\ref{thm:erst} and \ref{thm:hk}, the authors had the case $G_{\emptyset}=\{1\}$ in mind and, therefore, used a slightly different notion of curvature. But one can check that the theorems do also hold in our setting with arbitrary groups $G_{\emptyset}$ and Gersten-Stallings angles as defined in Section~\ref{sub:curvature}, see \cite{Cuno2011} for details.
\end{remark}
In the light of Theorem~\ref{thm:hk}, one may wonder about the Euclidean case. Let us therefore assume that the triangle of groups is Euclidean and ask under which conditions the colimit group $\mathcal{G}$ has a non-abelian free subgroup! A first class of examples to look at are \emph{Euclidean triangle groups}:
\[ \begin{array}{c} \Delta(k,l,m)\cong\left\langle a,b,c\,:\,a^{2},b^{2},c^{2},(ab)^{k},(ac)^{l},(bc)^{m}\right\rangle \smallskip \\ \text{with}~(k,l,m)\in\left\{(3,3,3),(2,4,4),(2,3,6)\right\} \end{array} \] 
Each of these groups happens to be the colimit group of the Euclidean triangle of groups based on the following data:
\[
\begin{array}{c} G_{\varnothing}\cong\{1\},\;G_{\{1\}}\cong\left\langle a\,:\,a^{2}=1\right\rangle,\;G_{\{2\}}\cong\left\langle b\,:\,b^{2}=1\right\rangle,\;G_{\{3\}}\cong\left\langle c\,:\,c^{2}=1\right\rangle, \smallskip \\ G_{\{1,2\}}\cong\left\langle a,b\,:\,a^{2}=b^{2}=(ab)^{k}=1\right\rangle,\;G_{\{1,3\}}\cong\left\langle a,c\,:\,a^{2}=c^{2}=(ac)^{l}=1\right\rangle, \smallskip \\ G_{\{2,3\}}\cong\left\langle b,c\,:\,b^{2}=c^{2}=(bc)^{m}=1\right\rangle \end{array}
\]
Here, as in Section~\ref{sub:example}, the homomorphisms $\varphi_{J_{1}J_{2}}:G_{J_{1}}\rightarrow G_{J_{2}}$ are implicitly given by $a\mapsto a$, $b\mapsto b$, $c\mapsto c$. It turns out that the Gersten-Stallings angles amount to $\sphericalangle_{\{1,2\}}=\nicefrac{\pi\,}{k}$, $\sphericalangle_{\{1,3\}}=\nicefrac{\pi\,}{l}$, $\sphericalangle_{\{2,3\}}=\nicefrac{\pi\,}{m}$, whence the triangle of groups is actually Euclidean. The algebraic structure of the colimit group $\mathcal{G}$ can be revealed by geometry.

\begin{figure}
\begin{center}
\input{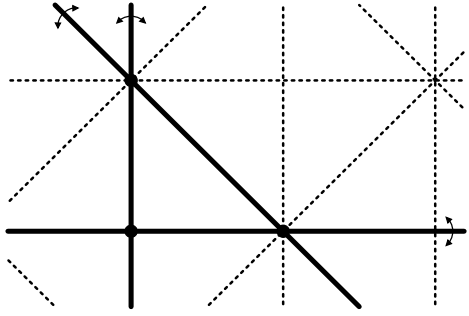}
\end{center}
\caption{The Euclidean triangle group $\Delta(2,4,4)$ realised as $S\leq\isom(\mathbb{E}^{2})$.}
\label{fig:trianglegroup}
\end{figure}

Consider three lines in the Euclidean plane $\mathbb{E}^{2}$ enclosing a triangle with angles $\nicefrac{\pi\,}{k}$, $\nicefrac{\pi\,}{l}$, $\nicefrac{\pi\,}{m}$, see Figure~\ref{fig:trianglegroup}. The reflections along these lines generate a subgroup $S\leq\isom(\mathbb{E}^{2})$ that is isomorphic to $\Delta(k,l,m)$. Since $\isom(\mathbb{E}^{2})$ is solvable, $S\leq\isom(\mathbb{E}^{2})$ must be solvable, too. Therefore, $\Delta(k,l,m)$ is solvable and cannot have a non-abelian free subgroup, see \cite{Cuno2011} for details.

In this section, we generalise the geometric approach. More precisely, we use a construction introduced by Bridson in \cite{Bridson1991} to study all \emph{non-degenerate} Euclidean triangles of groups, i.\,e.~all Euclidean triangles of groups with the property that each Gersten-Stallings angle is strictly between $0$ and $\pi$. Since the strict upper bound $\pi$ is our standing assumption, being non-degenerate actually means that none of the Gersten-Stallings angles is 0. Given such a triangle of groups, we will construct a simplicial complex $\mathcal{X}$. The action of $\mathcal{G}$ on $\mathcal{X}$ will give us new insight into the structure of $\mathcal{G}$. Actually, it turns out that if $\mathcal{X}$ \emph{branches}, i.\,e.~if $|\mathcal{X}|$ is not a topological manifold any more, then the colimit group $\mathcal{G}$ has a non-abelian free subgroup, see Theorem~\ref{thm:freesubgroup}.

This allows us to give an answer to a problem mentioned by Kopteva and Williams in \cite[p.\,58, l.\,24]{KoptevaWilliams2008}, who wondered if the class of colimit groups of non-spherical triangles of groups satisfies the Tits alternative.

As already mentioned, the construction and the basic properties of $\mathcal{X}$ have been introduced by Bridson in \cite{Bridson1991}. In Section~\ref{sub:bridson}, we summarise what is relevant for our work, and apply it from Section~\ref{sub:billiardstheorem} onwards. Our proofs are based on ideas and techniques that go back to two Diplomarbeiten under supervision of Bieri, namely by Lorenz in \cite{Lorenz1995} and Brendel in \cite{Brendel2004}. Both Lorenz and Brendel use altitudes in triangles to detect non-abelian free subgroups, but under additional assumptions on the Gersten-Stallings angles. We use the language of billiards instead, which gives us the flexibility we need.
\subsection{Bridson's simplicial complex} \label{sub:bridson}
Given a non-degenerate Euclidean triangle of groups, we define an abstract simplicial complex $\mathcal{X}$ as follows:
\[
\begin{array}{c@{\;}c@{\;}l}
0\text{-simplices} & := & \left\{\left\{g\mathcal{G}|_{\{1,2\}}\right\}\,:\,g\in\mathcal{G}\right\} \smallskip \\ && \;\sqcup\left\{\left\{g\mathcal{G}|_{\{1,3\}}\right\}\,:\,g\in\mathcal{G}\right\} \smallskip \\ && \;\sqcup\left\{\left\{g\mathcal{G}|_{\{2,3\}}\right\}\,:\,g\in\mathcal{G}\right\} \bigskip \\
1\text{-simplices} & := & \left\{\left\{g\mathcal{G}|_{\{1,2\}},g\mathcal{G}|_{\{1,3\}}\right\}\,:\,g\in\mathcal{G}\right\} \smallskip \\ && \;\sqcup\left\{\left\{g\mathcal{G}|_{\{1,2\}},g\mathcal{G}|_{\{2,3\}}\right\}\,:\,g\in\mathcal{G}\right\} \smallskip \\ && \;\sqcup\left\{\left\{g\mathcal{G}|_{\{1,3\}},g\mathcal{G}|_{\{2,3\}}\right\}\,:\,g\in\mathcal{G}\right\} \bigskip \\
2\text{-simplices} & := & \left\{\left\{g\mathcal{G}|_{\{1,2\}},g\mathcal{G}|_{\{1,3\}},g\mathcal{G}|_{\{2,3\}}\right\}\,:\, g\in\mathcal{G}\right\}
\end{array}
\]
\begin{remark} \label{rem:usefulmetric}
The same definition can, of course, be given for any triangle of groups. But we need it only for non-degenerate Euclidean triangles of groups.
\end{remark}
\subsubsection{Group action and stabilisers} \label{subsub:stabilisers}
\begin{figure}
\begin{center}
\input{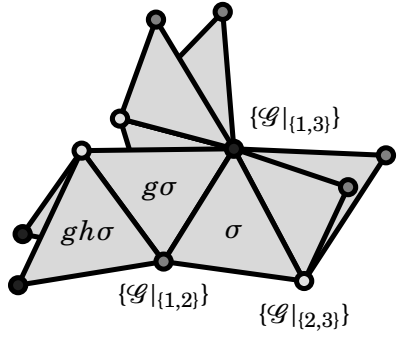}
\end{center}
\caption{A part of the simplicial complex $\mathcal{X}$ including translates of $\sigma$ by group elements $g$ and $gh$ with $g\in\mathcal{G}|_{\{1\}}$ and $h\in\mathcal{G}|_{\{2\}}$.}
\label{fig:complex}
\end{figure}

We will use the letter $\sigma$ to denote the 2-simplex that is represented by the identity, i.\,e.~$\sigma:=\{\mathcal{G}|_{\{1,2\}},\mathcal{G}|_{\{1,3\}},\mathcal{G}|_{\{2,3\}}\}$. There is a natural action of $\mathcal{G}$ on $\mathcal{X}$ given by left-multiplication of each coset. A fundamental domain for this action consists of the 2-simplex $\sigma$ and its faces, see Figure~\ref{fig:complex}.

Bridson mentioned in \cite[p.\,431, ll.\,8--9]{Bridson1991} that ``the pattern of stabilisers in this fundamental domain is precisely the original triangle of groups.'' Since we distinguish between the diagram and its image in the colimit group, we would replace ``original triangle of groups'' by ``pattern of subgroups $\mathcal{G}|_{K}$ with $K\subseteq\{1,2,3\}$ and $|K|\leq 2$.'' Anyway, notice that there is a very easy way to prove this observation using the intersection theorem. Indeed, we can observe that for all pairwise distinct elements $a,b,c\in\{1,2,3\}$:
\begin{enumerate}
\item $\begin{array}[t]{l} \operatorname{stab}_{\mathcal{G}}\left(\left\{\mathcal{G}|_{\{a,b\}}\right\}\right)=\left\{g\in\mathcal{G}\,:\,g\mathcal{G}|_{\{a,b\}}=\mathcal{G}|_{\{a,b\}}\right\}=\mathcal{G}|_{\{a,b\}} \end{array}$
\item $\begin{array}[t]{l} \operatorname{stab}_{\mathcal{G}}\left(\left\{\mathcal{G}|_{\{a,b\}},\mathcal{G}|_{\{a,c\}}\right\}\right) \smallskip \\ \;=\left\{g\in\mathcal{G}\,:\,g\mathcal{G}|_{\{a,b\}}=\mathcal{G}|_{\{a,b\}},\,g\mathcal{G}|_{\{a,c\}}=\mathcal{G}|_{\{a,c\}}\right\} \smallskip \\ \;=\mathcal{G}|_{\{a,b\}}\cap\mathcal{G}|_{\{a,c\}}=\mathcal{G}|_{\{a\}}\end{array}$
\item $\begin{array}[t]{l} \operatorname{stab}_{\mathcal{G}}\left(\left\{\mathcal{G}|_{\{a,b\}},\mathcal{G}|_{\{a,c\}},\mathcal{G}|_{\{b,c\}}\right\}\right) \smallskip \\ \;=\left\{g\in\mathcal{G}\,:\,g\mathcal{G}|_{\{a,b\}}=\mathcal{G}|_{\{a,b\}},\,g\mathcal{G}|_{\{a,c\}}=\mathcal{G}|_{\{a,c\}},\,g\mathcal{G}|_{\{b,c\}}=\mathcal{G}|_{\{b,c\}}\right\} \smallskip \\ \;=\mathcal{G}|_{\{a,b\}}\cap\mathcal{G}|_{\{a,c\}}\cap\mathcal{G}|_{\{b,c\}}=\mathcal{G}|_{\emptyset} \end{array}$
\end{enumerate}
Hence, by (1) to (3), the stabilisers of the $0$-simplices are the groups $\mathcal{G}|_{K}$ with $K\subseteq\{1,2,3\}$ and $|K|=2$. The stabilisers of the $1$-simplices and the $2$-simplex are their pairwise and triple intersections, which are precisely the groups $\mathcal{G}|_{K}$ with $K\subseteq\{1,2,3\}$ and $|K|=1$ and $|K|=0$.
\subsubsection{Simple connectedness} \label{subsub:simplyconnected}
Let us now consider the geometric realisation $|\mathcal{X}|$. As usual, it is equipped with the \emph{weak topology}. For details about abstract simplicial complexes and their geometric realisations we refer to \cite[Sections 1.1--1.3]{Munkres1984}.
\begin{lemma}[Behr]
\label{lem:behr}
The geometric realisation $|\mathcal{X}|$ is simply connected.
\end{lemma}
Behr has proved a slightly more general version of this lemma in 1975, see \cite[Satz 1.2]{Behr1975}. Roughly speaking, he translates products $h_{1}h_{2}\dotsb h_{n}$ of elements $h_{i}\in\mathcal{G}|_{K_{i}}$ with $K_{i}\subseteq\{1,2,3\}$ and $|K_{i}|=2$ to edge paths in $|\mathcal{X}|$, and vice versa.
\begin{remark}
In order to keep the notation simple, we use the same symbols to refer to simplices in $\mathcal{X}$ and to their geometric realisations in $|\mathcal{X}|$. Moreover, whenever we talk about a simplex in $|\mathcal{X}|$ without further specification, we mean the closed simplex.
\end{remark}
\subsubsection{Metric structure} \label{subsub:metric}
The geometric realisation $|\mathcal{X}|$ can be equipped with a piecewise Euclidean metric structure. The triangle of groups is non-degenerate and Euclidean. We may therefore pick a closed triangle $\Delta$ in the Euclidean plane $\mathbb{E}^{2}$ with the property that its angles agree with the Gersten-Stallings angles of the triangle of groups. For every subset $K\subseteq\{1,2,3\}$ with $|K|=2$ we label the corresponding vertex of $\Delta$ with angle $\sphericalangle_{K}$ by $K$. For later purposes, let us label the edges of $\Delta$, too. The edge between the two vertices that are labelled by $K_{1}$ and $K_{2}$ is labelled by their intersection $K_{1}\cap K_{2}$.
\begin{figure}
\begin{center}
\input{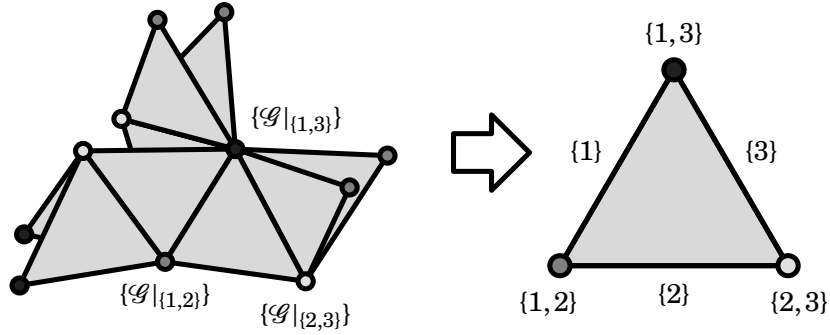}
\end{center}
\caption{Labelling the triangle $\Delta$ and constructing $h:|\mathcal{X}|\rightarrow\Delta$.}
\label{fig:sixparts}
\end{figure}

Next, we construct a continuous map $h:|\mathcal{X}|\rightarrow\Delta$. Every $0$-simplex is of the form $\{g\mathcal{G}|_{K}\}$ for some $K\subseteq\{1,2,3\}$ with $|K|=2$. Map it to the vertex that is labelled by $K$. In order to continue this map to the higher dimensional simplices, map the $1$-simplices homeomorphically to the \emph{corresponding} edges, i.\,e.~if a $1$-simplex is of the form $\{g\mathcal{G}|_{\{a,b\}},g\mathcal{G}|_{\{a,c\}}\}$, map it to the edge that is labelled by $\{a,b\}\cap\{a,c\}=\{a\}$. For every $2$-simplex $\tau$, use Schoenflies' Theorem to continue the homeomorphism $h|_{\partial\tau}:\partial\tau\rightarrow\partial\Delta$ to $h|_{\tau}:\tau\rightarrow\Delta$. The latter ones assemble to the desired continuous map $h:|\mathcal{X}|\rightarrow\Delta$. Given two points $x,y\in\tau$ we can now measure their \emph{local distance}:
\[ \dist_{\tau}(x,y):=|\hspace{-1.5pt}|\,h|_{\tau}(x)-h|_{\tau}(y)\,|\hspace{-1.5pt}| \]
The local distance $\dist_{\tau}:\tau\times\tau\rightarrow\mathbb{R}$ is a metric on the single $2$-simplex $\tau$. The geometric realisation $|\mathcal{X}|$ equipped with all these local distances is called an \emph{$\mathbb{E}$-simplicial complex}, see \cite[Section 1.1]{Bridson1991} for a formal definition. In order to extend them to a metric on $|\mathcal{X}|$, we follow Bridson's work.
\begin{definition}[``$m$-chain'']
Let $x,y\in|\mathcal{X}|$. An $m$-chain from $x$ to $y$ is a finite sequence $\mathcal{C}=(x_{0},x_{1},\dotsc,x_{m})$ of points in $|\mathcal{X}|$ with the property that $x_{0}=x$ and $x_{m}=y$, and that for every $1\leq i\leq m$ both $x_{i-1}$ and $x_{i}$ are contained in some common $2$-simplex $\tau_{i}$.
\end{definition}
Let $\mathcal{C}$ be an $m$-chain as above. Then, the \emph{length} of $\mathcal{C}$ is defined by:
\[ \operatorname{length}(\mathcal{C}):=\sum_{i=1}^{m}\dist_{\tau_{i}}(x_{i-1},x_{i}) \]
For every $1\leq i\leq m$ there is a unique geodesic from $x_{i-1}$ to $x_{i}$ in $\tau_{i}$. We call it the \emph{segment} from $x_{i-1}$ to $x_{i}$. The concatenation of all segments is called the \emph{path} induced by $\mathcal{C}$. It is denoted by $\pathbrackets{\mathcal{C}}$.
\begin{remark}
Notice that neither $\operatorname{length}(\lambda)$ nor $\pathbrackets{\mathcal{C}}$ depends on the choice of $\tau_{i}$, i.\,e.~if there are two $2$-simplices $\tau_{i}$ and $\widetilde{\tau}_{i}$ such that $x_{i-1},x_{i}\in\tau_{i}\cap\widetilde{\tau}_{i}$, then the local distances $\dist_{\tau_{i}}(x_{i-1},x_{i})$ and $\dist_{\widetilde{\tau}_{i}}(x_{i-1},x_{i})$ agree and the segment $\pathbrackets{x_{i-1},x_{i}}$ is well defined.
\end{remark}
Let $x,y\in|\mathcal{X}|$. Since $|\mathcal{X}|$ is path connected, there is a path from $x$ to $y$. By the construction given in Behr's proof or by a direct argument, we can even find a path from $x$ to $y$ that is induced by some $m$-chain $\mathcal{C}$. Hence, there is a function $\dist:|\mathcal{X}|\times|\mathcal{X}|\rightarrow\mathbb{R}$ given by:
\[ \dist(x,y):=\inf\left\{\operatorname{length}(\mathcal{C})\,:\,\text{$\exists\,m$ such that $\mathcal{C}$ is an $m$-chain from $x$ to $y$}\right\} \]
It is straightforward to see that $\dist:|\mathcal{X}|\times|\mathcal{X}|\rightarrow\mathbb{R}$ is a pseudometric. On the other hand, \emph{distinguishability}, i.\,e.~$\dist(x,y)=0$ implies $x=y$, is an issue. At this point, a lemma from \cite[Section 1.2]{Bridson1991} comes into play. It uses that $|\mathcal{X}|$ is an $\mathbb{E}$-simplicial complex with a finite set of shapes, i.\,e.~with a finite set of isometry classes of simplices.
\begin{lemma}[Bridson] \label{lem:localagreement}
\wegschummeln{9}For every $x\in|\mathcal{X}|$ there is an $\varepsilon(x)>0$ with the following property: For every $y\in|\mathcal{X}|$ with $\dist(x,y)<\varepsilon(x)$ there is a common $2$-simplex $\tau$ containing both $x$ and $y$ such that the distances $\dist_{\tau}(x,y)$ and $\dist(x,y)$ agree. \hfill $\Box$
\end{lemma}
This lemma immediately implies distinguishabilty, whence the pseudometric $\dist:|\mathcal{X}|\times|\mathcal{X}|\rightarrow\mathbb{R}$ is actually a metric. Even more, it makes $|\mathcal{X}|$ a complete geodesic metric space. For details, see \cite[Theorem 1.1]{Bridson1991}.

\begin{remark}
As mentioned in \cite[p.\,381, ll.\,25--31]{Bridson1991}, the topology induced by the metric $\dist:|\mathcal{X}|\times|\mathcal{X}|\rightarrow\mathbb{R}$ is coarser, and may even be strictly coarser, than the weak topology. However, $|\mathcal{X}|$ remains simply connected as a metric space.
\end{remark}

\begin{remark}[``arc length''] \label{rem:lengthsagree}
Another important application of Lemma~\ref{lem:localagreement} concerns the arc length of paths. Given an $m$-chain $\mathcal{C}=(x_{0},x_{1},\dotsc,x_{m})$, the metric $\dist:|\mathcal{X}|\times|\mathcal{X}|\rightarrow\mathbb{R}$ allows us to determine the arc length of the path $\pathbrackets{\mathcal{C}}$, see e.\,g.~\cite[Definition I.1.18]{BridsonHaefliger1999}. In the light of Lemma~\ref{lem:localagreement}, it is easy to verify that the arc length of the path $\pathbrackets{\mathcal{C}}$ agrees with $\operatorname{length}(\mathcal{C})$.
\end{remark}
\subsubsection{CAT(0) property} \label{subsub:catzero}
From now on, we will consider $|\mathcal{X}|$ as a metric space. A crucial observation is that $|\mathcal{X}|$ has the $\operatorname{CAT}(0)$ property. In order to prove this, we will verify the \emph{link condition}.
\begin{definition}[``geometric link'']
\label{def:link}
Let $x\in|\mathcal{X}|$. The (geometric) closed star $\operatorname{St}(x)$ is the union of the geometric relisations of all simplices that contain $x$. If $y\in\operatorname{St}(x)\setminus\{x\}$, then there is at least one $2$-simplex that contains both $x$ and $y$. So, we may consider the segment $\pathbrackets{x,y}$. Two such segments $\pathbrackets{x,y}$ and $\pathbrackets{x,y'}$ are called equivalent if one of them is contained in the other. We call the set of equivalence classes the geometric link of $x$. It is denoted by $\operatorname{Lk}(x,|\mathcal{X}|)$.
\end{definition}
The geometric link $\operatorname{Lk}(x,|\mathcal{X}|)$ can be equipped with a metric structure. First, we consider certain subsets of $\operatorname{Lk}(x,|\mathcal{X}|)$. For every $2$-simplex $\tau$ in $\operatorname{St}(x)$ let $\operatorname{Lk}(x,\tau)$ be the subset of all elements of $\operatorname{Lk}(x,|\mathcal{X}|)$ that are represented by segments in $\tau$. Notice that, as soon as one representative has this property, all representatives do. The subset $\operatorname{Lk}(x,\tau)\subseteq\operatorname{Lk}(x,|\mathcal{X}|)$ has a natural metric $\dist_{\operatorname{Lk}(x,\tau)}:\operatorname{Lk}(x,\tau)\times\operatorname{Lk}(x,\tau)\rightarrow\mathbb{R}$ given by the Euclidean angle:
\[ \dist_{\operatorname{Lk}(x,\tau)}\left(\pathbrackets{x,y}_{\sim},\pathbrackets{x,y'}_{\sim}\right):=\angle_{h|_{\tau}(x)}\left(h|_{\tau}(y),h|_{\tau}(y')\right)\in[0,\pi] \]
If $x$ is in the $1$-skeleton of $|\mathcal{X}|$, then every $(\operatorname{Lk}(x,\tau),\dist_{\operatorname{Lk}(x,\tau)})$ is isometrically isomorphic to a closed interval of length $\sphericalangle_{\{1,2\}}$, $\sphericalangle_{\{1,3\}}$, $\sphericalangle_{\{2,3\}}$, or $\pi$, see \circled{1} and \circled{2} in Figure~\ref{fig:links}. In particular, we may interpret the subsets $\operatorname{Lk}(x,\tau)\subseteq\operatorname{Lk}(x,|\mathcal{X}|)$ as $1$-simplices and, at least after a barycentric subdivision of each simplex, the whole geometric link $\operatorname{Lk}(x,|\mathcal{X}|)$ as a simplicial complex.

Even more, it is an $\mathbb{E}$-simplicial complex with a finite set of shapes. As in Section~\ref{subsub:metric}, the connected components of $\operatorname{Lk}(x,|\mathcal{X}|)$ can be equipped with a pseudometric. This pseudometric turns out to be a metric which makes every connected component a complete geodesic metric space. If we set the distance of elements from distinct connected components to $\infty$, we obtain an extended metric $\dist_{\operatorname{Lk}(x,|\mathcal{X}|)}:\operatorname{Lk}(x,|\mathcal{X}|)\times\operatorname{Lk}(x,|\mathcal{X}|)\rightarrow\mathbb{R}\cup\{\infty\}$.

If $x$ is not in the $1$-skeleton of $|\mathcal{X}|$, then it must be in the interior of some $2$-simplex $\tau$. In this case, $(\operatorname{Lk}(x,\tau),\dist_{\operatorname{Lk}(x,\tau)})$ is isometrically isomorphic to the standard $1$-sphere $\mathbb{S}^{1}$, see \circled{3} in Figure~\ref{fig:links}. Since $\operatorname{Lk}(x,\tau)=\operatorname{Lk}(x,|\mathcal{X}|)$, the metric $\dist_{\operatorname{Lk}(x,\tau)}$ is actually a metric $\dist_{\operatorname{Lk}(x,|\mathcal{X}|)}:\operatorname{Lk}(x,|\mathcal{X}|)\times\operatorname{Lk}(x,|\mathcal{X}|)\rightarrow\mathbb{R}$.

For more details, in particular for a remark about equivalent definitions, we refer to \cite[Sections 1.7.14--1.7.15]{BridsonHaefliger1999}. Notice that Bridson and Haefliger consider the open star instead on the closed one. But, in the end, this does not make a difference.
\begin{figure}
\begin{center}
\input{img-16.pspdftex}
\end{center}
\caption{The sets $\operatorname{Lk}(x,\tau)$ for different points $x\in\tau$.}
\label{fig:links}
\end{figure}
\begin{definition}[``link condition'']\schummeln{8}The geometric realisation $|\mathcal{X}|$ satisfies the link condition if for every $x\in|\mathcal{X}|$ and every pair of points $a,b\in\operatorname{Lk}(x,|\mathcal{X}|)$ with $\dist_{\operatorname{Lk}(x,|\mathcal{X}|)}(a,b)<\pi$ there is a unique geodesic from $a$ to $b$. Or, equivalently, if every injective loop $\lambda:\mathbb{S}^{1}\rightarrow\operatorname{Lk}(x,|\mathcal{X}|)$ has arc length at least $2\pi$.
\end{definition}
The equivalence relies on the fact that $|\mathcal{X}|$ is a $2$-dimensional $\mathbb{E}$-simplicial complex with a finite set of shapes. One can either prove it directly or apply \cite[Theorem I.7.55, (3)\,$\Leftrightarrow$\,(1)]{BridsonHaefliger1999}, \cite[Theorem II.5.5, (3)\,$\Leftrightarrow$\,(2)]{BridsonHaefliger1999}, \cite[Lemma II.5.6]{BridsonHaefliger1999}.

\begin{lemma}[Bridson and Gersten-Stallings]
The geometric realisation $|\mathcal{X}|$ satisfies the link condition. 
\end{lemma}
A proof of this lemma has been given by Bridson in \cite[p.\,431, ll.\,19--25]{Bridson1991} and by Gersten and Stallings in \cite[p.\,499, ll.\,19--28]{Stallings1991}. The argument is both simple and important, so let us outline the main ideas! \bigskip \\
\textbf{Proof sketch.} The most interesting case occurs, when $x$ is a $0$-simplex of $|\mathcal{X}|$, w.\,l.\,o.\,g.~$x=\{g\mathcal{G}|_{\{1,2\}}\}$. Then, every injective loop $\lambda:\mathbb{S}^{1}\rightarrow\operatorname{Lk}(x,|\mathcal{X}|)$ traverses some finite number of $1$-simplices, say $m$. If the first one is $\operatorname{Lk}(x,g\sigma)$, the next ones are $\operatorname{Lk}(x,gh_{1}h_{2}\dotsb h_{i}\sigma)$ with elements $h_{i}$ alternately in $\mathcal{G}|_{\{1\}}\setminus\mathcal{G}|_{\emptyset}$ and $\mathcal{G}|_{\{2\}}\setminus\mathcal{G}|_{\emptyset}$. At the end, $\lambda$ traverses $\operatorname{Lk}(x,gh_{1}h_{2}\dotsb h_{m-1}\sigma)$. Since $\lambda$ is a loop, we know that there is an $h_{m}$ such that $\operatorname{Lk}(x,gh_{1}h_{2}\dotsb h_{m}\sigma)=\operatorname{Lk}(x,g\sigma)$, which is equivalent to $h_{1}h_{2}\dotsb h_{m}\in\mathcal{G}|_{\emptyset}$. Let us write $h:=h_{1}h_{2}\dotsb h_{m}$. We may assume w.\,l.\,o.\,g.~that $h_{1},h_{3},\dotsc,h_{m-1}\in\mathcal{G}|_{\{1\}}\setminus\mathcal{G}|_{\emptyset}$ and $h_{2},h_{4},\dotsc,h_{m}\in\mathcal{G}|_{\{2\}}\setminus\mathcal{G}|_{\emptyset}$. Since $h\in\mathcal{G}|_{\emptyset}$, also $h_{m}h^{-1}\in\mathcal{G}|_{\{2\}}\setminus\mathcal{G}|_{\emptyset}$. Now, we proceed similarly to (\ref{item:swap}) in the proof of the intersection theorem. The preimages:
\[ \begin{array}{l} {\nu_{\{1\}}}^{-1}(h_{1}),{\nu_{\{1\}}}^{-1}(h_{3}),\dotsc,{\nu_{\{1\}}}^{-1}(h_{m-1})\in G_{\{1\}}\setminus\varphi_{\emptyset\{1\}}(G_{\emptyset}) \smallskip \\ {\nu_{\{2\}}}^{-1}(h_{2}),{\nu_{\{2\}}}^{-1}(h_{4}),\dotsc,{\nu_{\{2\}}}^{-1}(h_{m}h^{-1})\in G_{\{2\}}\setminus\varphi_{\emptyset\{2\}}(G_{\emptyset}) \end{array} \]
assemble to an element ${\nu_{\{1\}}}^{-1}(h_{1})\cdot{\nu_{\{2\}}}^{-1}(h_{2})\cdot\dotsc\cdot{\nu_{\{2\}}}^{-1}(h_{m}h^{-1})$ of the amalgamated free product $G_{\{i\}}\ast_{G_{\varnothing}}G_{\{j\}}$. It is in the kernel of the homomorphism $\alpha:G_{\{i\}}\ast_{G_{\varnothing}}G_{\{j\}}\rightarrow G_{\{i,j\}}$ introduced in Section~\ref{sub:curvature}. But, by the normal form theorem, see \cite[Lemma 1]{Miller1968}, it is non-trivial. Therefore, $m\geq\nicefrac{2\pi\,}{\sphericalangle_{\{1,2\}}}$.

Remark~\ref{rem:lengthsagree}, which does also hold for $m$-chains in $\operatorname{Lk}(x,|\mathcal{X}|)$, implies that the arc length of $\lambda$ is equal to $m\cdot\sphericalangle_{\{1,2\}}$, which can be estimated from below as follows  $m\cdot\sphericalangle_{\{1,2\}}\geq\nicefrac{2\pi\,}{\sphericalangle_{\{1,2\}}}\cdot\sphericalangle_{\{1,2\}}=2\pi$, whence we are done. The link condition for the other cases, i.\,e.~if $x$ is in the interior of a $1$-simplex or in the interior of a $2$-simplex, is almost immediate. \hfill ``$\Box$''\bigskip \\
Once we have convinced ourselves that the geometric realisation $|\mathcal{X}|$ satisfies the link condition, we may apply Bridson's main theorem \cite[Section 2, Main Theorem, (11)\,$\Rightarrow$\,(2)]{Bridson1991}.
\begin{theorem}[Bridson]
Since $|\mathcal{X}|$ is a simply connected $\mathbb{E}$-simplicial complex with a finite set of shapes that satisfies the link condition, it has the $\operatorname{CAT}(0)$ property. \hfill $\Box$
\end{theorem}
\subsubsection{Geodesics} \label{subsub:geodesics}
Let $\pathbrackets{\mathcal{C}}$ be the path induced by an $m$-chain $\mathcal{C}=(x_{0},x_{2},\dotsc,x_{m})$. There is a necessary and, as we will see in Lemma~\ref{lem:localtoglobal}, sufficient condition for $\pathbrackets{\mathcal{C}}$ to be a geodesic, namely that $\mathcal{C}$ is \emph{straight}, i.\,e.~that there are no obvious shortcuts at the points $x_{1},x_{3},\dotsc,x_{m-1}$. Let us make the notion of straightness a little more precise!
\begin{definition}[``straight $m$-chain'']\schummeln{4}
An $m$-chain $\mathcal{C}=(x_{0},x_{1},\dotsc,x_{m})$ is called straight if for every $1\leq i\leq m-1$ the distance between $\pathbrackets{x_{i},x_{i-1}}_{\sim}$ and $\pathbrackets{x_{i},x_{i+1}}_{\sim}$ in $\operatorname{Lk}(x_{i},|\mathcal{X}|)$ is at least $\pi$.
\end{definition}
\begin{remark} \label{rem:specialcase}
In general, at least when $x_{i}$ is a $0$-simplex, it is not easy to determine the distance between $\pathbrackets{x_{i},x_{i-1}}_{\sim}$ and $\pathbrackets{x_{i},x_{i+1}}_{\sim}$ in $\operatorname{Lk}(x_{i},|\mathcal{X}|)$. But the link condition ensures that, once we are able to connect them by an injective path of arc length $\pi$, the distance between them is actually equal to $\pi$.
\end{remark}
At this point, the $\operatorname{CAT}(0)$ property comes into play. It allows us to conclude from the local property ``straight $m$-chain'' to the global property ``geodesic''.
\begin{lemma} \label{lem:localtoglobal}
If $\mathcal{C}$ is a straight $m$-chain, then $\pathbrackets{\mathcal{C}}$ is a geodesic.
\end{lemma}
\textbf{Proof sketch.} A proof can be given by showing that every straight $m$-chain induces a local geodesic, which is an easy consequence of \cite[Section 2, Main Theorem, (2)\,$\Rightarrow$\,(5)]{Bridson1991}. Finally, \cite[Proposition II.1.4 (2)]{BridsonHaefliger1999} tells us that every local geodesic is a geodesic. \hfill ``$\Box$''\bigskip\\
Moreover, the geodesics in $|\mathcal{X}|$ are unique, see \cite[Section 2, Main Theorem, (2)\,$\Rightarrow$\,(1)]{Bridson1991}. This will be of relevance in the proof of the billiards theorem.
\subsection{Billiards theorem} \label{sub:billiardstheorem}
In this section, we still assume given a non-degerate Euclidean triangle of groups and consider \emph{billiard shots} and \emph{billiard sequences} on the triangle $\Delta$.
\subsubsection{Billiard shots and billiard sequences}
\begin{figure}
\begin{center}
\input{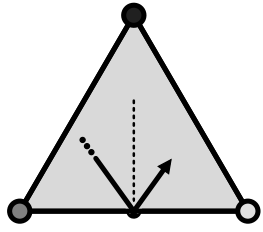}
\end{center}
\caption{The reflection rule for billiard shots.}
\label{fig:backandforth}
\end{figure}

We choose some point $y_{0}$ in the interior of $\Delta$ and some direction. Then, we consider the path that starts at $y_{0}$ and goes in a straight line into the chosen direction. Eventually, this path leaves $\Delta$. Let $y_{1}\in\partial\Delta$ be its last point in $\Delta$. If this point is a vertex (``the ball is in the pocket''), we withdraw the path. Otherwise, it is in the interior of an edge (``the ball hits the cushion''), which allows us to reflect the path according to the rule that the \emph{angle of incidence} is equal to the \emph{angle of reflection}, see Figure~\ref{fig:backandforth}. Now, we can go on. Whenever the path leaves $\Delta$ at some point in the interior of an edge, we reflect it again. After some finite number of reflections, say at the points $y_{1},y_{2},\dotsc,y_{m-1}\in\partial\Delta$, we stop at some point $y_{m}$ in the interior of $\Delta$. The sequence $\mathcal{B}=(y_{0},y_{1},\dotsc,y_{m})$ is called a \emph{billiard sequence}, the induced path $\pathbrackets{\mathcal{B}}=\pathbrackets{y_{1},y_{2},\dotsc,y_{m}}$ is called a \emph{billiard shot}.
\begin{definition}[``adapted''] \label{def:adapted}
Given a billiard sequence $\mathcal{B}=(y_{0},y_{1},\dotsc,y_{m})$, we call an element $g\in\mathcal{G}$ adapted to $\mathcal{B}$, if it is a product $g_{1}g_{2}\dotsb g_{m-1}$ such that each factor $g_{i}\in\mathcal{G}|_{\{a_{i}\}}\setminus\mathcal{G}|_{\emptyset}$, where $\{a_{i}\}$ is the label of the edge whose interior contains $y_{i}$.
\end{definition}
\subsubsection{Statement and proof of the billiards theorem} \label{subsub:billiards}
\begin{theorem}
\label{thm:billiard}
Assume given a non-degenerate Euclidean triangle of groups and a closed triangle $\Delta$ in the Euclidean plane $\mathbb{E}^{2}$ with the property that its angles agree with the Gersten-Stallings angles of the triangle of groups. If an element $g\in\mathcal{G}$ is adapted to a billiard sequence $\mathcal{B}=(y_{0},y_{1},\dotsc,y_{m})$ on $\Delta$ with at least one reflection, i.\,e.~with $m\geq 2$, then $g$ is non-trivial.
\end{theorem}
\textbf{Proof.} The idea of the proof is to lift the billiard shot $\pathbrackets{\mathcal{B}}$ to the geometric realisation $|\mathcal{X}|$. For every $1\leq i\leq m$ we use $h|_{g_{1}g_{2}\dotsb g_{i-1}\sigma}:g_{1}g_{2}\dotsb g_{i-1}\sigma\rightarrow\Delta$ to lift the segment $\pathbrackets{y_{i-1},y_{i}}$. Let us make some observations!

\begin{enumerate}
\item \textbf{These lifts assemble to a path in $|\mathcal{X}|$.} For every $1\leq i\leq m-1$ the first segment $\pathbrackets{y_{i-1},y_{i}}$ is lifted by $h|_{g_{1}g_{2}\dotsb g_{i-1}\sigma}$, the second segment $\pathbrackets{y_{i},y_{i+1}}$ is lifted by $h|_{g_{1}g_{2}\dotsb g_{i}\sigma}$. To show that these two lifts actually fit together, we convince ourselves that in either case $y_{i}$ is lifted to the same point.
If the edge whose interior contains $y_{i}$ is labelled by $\{a\}$, then the two adjacent vertices are labelled by $\{a,b\}$ and $\{a,c\}$, where $b$ and $c$ are the remaining two elements of $\{1,2,3\}$. We can observe: 
\[ 
\begin{array}{l}
\left(h|_{g_{1}g_{2}\dotsb g_{i-1}\sigma}\right)^{-1}(y_{i})\in\left\{g_{1}g_{2}\dotsb g_{i-1}\mathcal{G}|_{\{a,b\}},g_{1}g_{2}\dotsb g_{i-1}\mathcal{G}|_{\{a,c\}}\right\} \smallskip \\
\left(h|_{g_{1}g_{2}\dotsb g_{i}\sigma}\right)^{-1}(y_{i})\in\left\{g_{1}g_{2}\dotsb g_{i}\mathcal{G}|_{\{a,b\}},g_{1}g_{2}\dotsb g_{i}\mathcal{G}|_{\{a,c\}}\right\}
\end{array}
\]
Since $g_{i}\in\mathcal{G}|_{\{a\}}$, the two $1$-simplices agree. Call them $\tau$ and observe:
\[
\left(h|_{g_{1}g_{2}\dotsb g_{i-1}\sigma}\right)^{-1}(y_{i})=\big(h|_{\tau}\big)^{-1}(y_{i})=\left(h|_{g_{1}g_{2}\dotsb g_{i}\sigma}\right)^{-1}(y_{i})
\]
So, the lift of $y_{i}$ is well defined. Let us denote it by $x_{i}$. For the lifts of the extremal points $y_{0}$ and $y_{m}$ we define analogously:
\[ 
x_{0}:=\big(h|_{\sigma}\big)^{-1}(y_{0})\quad\text{and}\quad x_{m}:=\left(h|_{g\sigma}\right)^{-1}(y_{m})
\]
By the above, we know that the lifts of the segments $\pathbrackets{y_{i-1},y_{i}}$ and $\pathbrackets{y_{i},y_{i+1}}$ assemble to the path $\pathbrackets{x_{i-1},x_{i},x_{i+1}}$ and, more general, that the lifts of all segments $\pathbrackets{y_{0},y_{1}},\pathbrackets{y_{1},y_{2}},\dotsc,\pathbrackets{y_{m-1},y_{m}}$ assemble to the path $\pathbrackets{\mathcal{C}}$ induced by the $m$-chain $\mathcal{C}:=\left(x_{0},x_{1},\dotsc,x_{m}\right)$.~\Checkmark
\item \textbf{The $m$-chain $\mathcal{C}$ is straight.} Let $1\leq i\leq m-1$. We construct a $2$-chain $\mathcal{L}$ from $\pathbrackets{x_{i},x_{i-1}}_{\sim}$ to $\pathbrackets{x_{i},x_{i+1}}_{\sim}$ in $\operatorname{Lk}(x_{i},|\mathcal{X}|)$ such that $\operatorname{length}(\mathcal{L})=\pi$ and the path $\pathbrackets{\mathcal{L}}$ is injective. Then, by Remark~\ref{rem:lengthsagree}, the path $\pathbrackets{\mathcal{L}}$ is of arc length $\pi$ and, by Remark~\ref{rem:specialcase}, the distance between $\pathbrackets{x_{i},x_{i-1}}_{\sim}$ and $\pathbrackets{x_{i},x_{i+1}}_{\sim}$ in $\operatorname{Lk}(x_{i},|\mathcal{X}|)$ is equal to $\pi$.

First, consider the edge whose interior contains $y_{i}$ and choose another point $\widetilde{y}_{i}$ in the interior of this edge. Let $\widetilde{x}_{i}$ be its lift:
\[
\widetilde{x}_{i}:=\left(h|_{g_{1}g_{2}\dotsb g_{i-1}\sigma}\right)^{-1}(\widetilde{y}_{i})=\left(h|_{g_{1}g_{2}\dotsb g_{i}\sigma}\right)^{-1}(\widetilde{y}_{i})
\]
Then, move to the geometric link $\operatorname{Lk}(x_{i},|\mathcal{X}|)$ and construct the $2$-chain $\mathcal{L}:=(\pathbrackets{x_{i},x_{i-1}}_{\sim},\pathbrackets{x_{i},\widetilde{x}_{i}}_{\sim},\pathbrackets{x_{i},x_{i+1}}_{\sim})$, see Figure~\ref{fig:liftingalpha}.
\begin{figure}
\begin{center}
\input{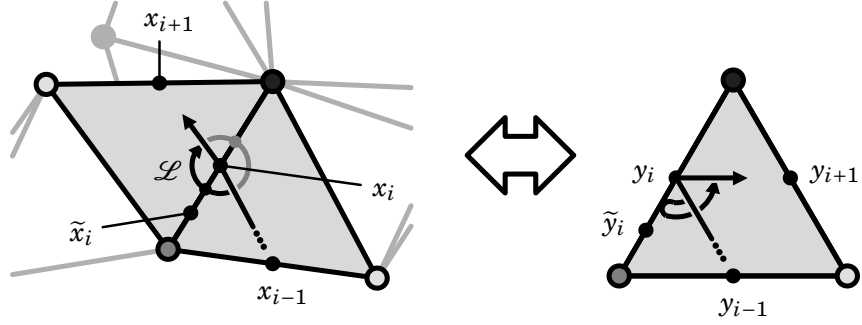}
\end{center}
\caption{Lift the billiard shot to $g_{1}g_{2}\dotsb g_{i-1}\sigma$ and $g_{1}g_{2}\dotsb g_{i}\sigma$.}
\label{fig:liftingalpha}
\end{figure}

Observe that the path $\pathbrackets{\mathcal{L}}$ traverses the interval $\operatorname{Lk}(x_{i},g_{1}g_{2}\dotsb g_{i-1}\sigma)$ until it reaches its endpoint $\pathbrackets{x_{i},\widetilde{x}_{i}}_{\sim}$. Then, it traverses the interval $\operatorname{Lk}(x_{i},g_{1}g_{2}\dotsb g_{i}\sigma)$. Therefore:
\[
\begin{array}{r@{\;}c@{\;}l}
\operatorname{length}(\mathcal{L}) & = & \dist_{\operatorname{Lk}(x_{i},g_{1}g_{2}\dotsb g_{i-1}\sigma)}\left(\pathbrackets{x_{i},x_{i-1}}_{\sim},\pathbrackets{x_{i},\widetilde{x}_{i}}_{\sim}\right) \smallskip \\
&& \,+\,\dist_{\operatorname{Lk}(x_{i},g_{1}g_{2}\dotsb g_{i}\sigma)}\left(\pathbrackets{x_{i},\widetilde{x}_{i}}_{\sim},\pathbrackets{x_{i},x_{i+1}}_{\sim}\right) \medskip \\
& = & \angle_{y_{i}}\left(y_{i-1},\widetilde{y}_{i}\right)+\angle_{y_{i}}\left(\widetilde{y}_{i},y_{i+1}\right)=\pi
\end{array}
\]
Since $g_{i}\not\in\mathcal{G}|_{\emptyset}$, the two intervals traversed by $\pathbrackets{\mathcal{L}}$ are actually not the same. So, $\pathbrackets{\mathcal{L}}$ must be injective.~\Checkmark
\end{enumerate}
With these two observations in mind the final conclusion that $g$ is non-trivial in $\mathcal{G}$ is almost immediate. By Lemma~\ref{lem:localtoglobal}, $\pathbrackets{\mathcal{C}}$ is a geodesic from $x_{0}\in\sigma$ to $x_{m}\in g\sigma$. Before going on, observe that any two points in $\sigma$ can be connected by a $1$-chain, which is, of course, straight and therefore induces a geodesic. But, as mentioned in Section~\ref{subsub:geodesics}, geodesics are unique. Hence, the unique geodesic between any two points in $\sigma$ is completely contained in $\sigma$. Let us now go back to our situation! We assume that $m\geq 2$. So, the geodesic $\pathbrackets{\mathcal{C}}$ leaves $\sigma$ eventually and, therefore, does not end in $\sigma$, i.\,e.~$x_{m}\not\in\sigma$, which implies that $g\sigma\neq\sigma$ and, finally, $g\neq 1$. \hfill $\Box$
\subsubsection{A first example} \label{subsub:playingbilliardszero}
We conclude with an example that illustrates the application of the billiards theorem. Assume given a triangle of groups with Gersten-Stallings angles $\sphericalangle_{\{1,2\}}=\sphericalangle_{\{1,3\}}=\sphericalangle_{\{2,3\}}=\nicefrac{\pi\,}{3}$. Since none of them is equal to $0$, it is easy to see that for every $a\in\{1,2,3\}$ there is an element $g_{a}\in\mathcal{G}|_{\{a\}}\setminus\mathcal{G}|_{\emptyset}$.

\begin{figure}
\begin{center}
\input{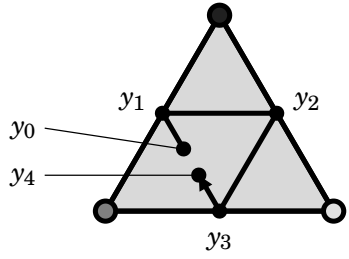}
\end{center}
\caption{A first billiard shot.}
\label{fig:firstbilliard}
\end{figure}

Their product $g:=g_{1}g_{2}g_{3}\in\mathcal{G}$, which has been studied in Theorem~\ref{thm:erst}, is adapted to the billiard sequence $\mathcal{B}=(y_{0},y_{1},y_{2},y_{3},y_{4})$ drawn in Figure~\ref{fig:firstbilliard} and is therefore non-trivial. Even more, we may continue the billiard shot. This yields billiard sequences of the form $(y_{0},y_{1},y_{2},y_{3},y_{1},y_{2},y_{3},\dotsc,y_{1},y_{2},y_{3},y_{4})$. Every \emph{power} of $g$, i.\,e.~every element $g^{n}\in\mathcal{G}$ with $n\in\mathbb{N}$, is adapted to such a billiard sequence and is therefore non-trivial. Hence, $g$ has infinite order.
\subsection{Constructing non-abelian free subgroups} \label{sub:playingbilliards}
Assume given a non-degenerate Euclidean triangle of groups. In particular, as in the example, for every $a\in\{1,2,3\}$ there is an element $g_{a}\in\mathcal{G}|_{\{a\}}\setminus\mathcal{G}|_{\emptyset}$. Let us recall the notion of branching! We say that the simplicial complex $\mathcal{X}$ \emph{branches} if the geometric realisation $|\mathcal{X}|$ is not a topological manifold any more. It is easy to see that $\mathcal{X}$ branches if and only if there is an $a\in\{1,2,3\}$ such that the index of $\mathcal{G}|_{\emptyset}$ in $\mathcal{G}|_{\{a\}}$ is at least 3 or there are two distinct $a,b\in\{1,2,3\}$ such that $\mathcal{G}|_{\{a,b\}}$ is not generated by $\mathcal{G}|_{\{a\}}$ and $\mathcal{G}|_{\{b\}}$.

\begin{figure}
\begin{center}
\input{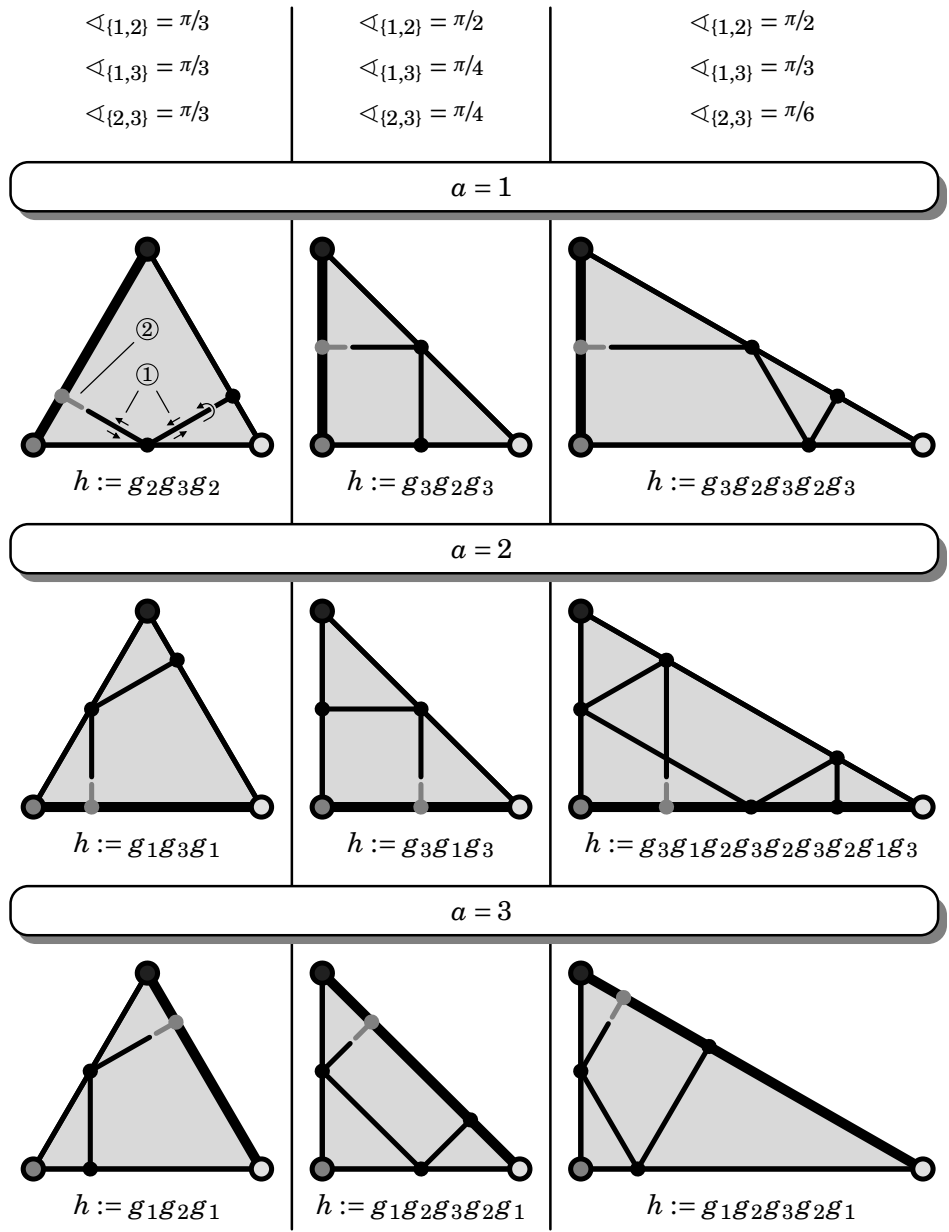}
\end{center}
\caption{The nine cases.}
\label{fig:ninecases}
\end{figure}

The following theorem says that branching already implies the existence of non-abelian free subgroups in $\mathcal{G}$.
\begin{theorem} \label{thm:freesubgroup}
Assume given a non-degenerate Euclidean triangle of groups. If the simplicial complex $\mathcal{X}$ branches, the colimit group $\mathcal{G}$ has a non-abelian free subgroup.
\end{theorem}
\begin{remark}
\label{rem:bassserre}
Notice that Theorem~\ref{thm:freesubgroup} is the $2$-dimensional analogue of a well known fact. Consider an amalgamated free product $X\ast_{A}Y$ with the property that the image of $A$ in $X$ and the image of $A$ in $Y$ have index at least 2. The associated Bass-Serre tree $\mathcal{T}$ branches if and only if one of the indices is at least $3$. In this case, the amalgamated free product $X\ast_{A}Y$ has a non-abelian free subgroup.
\end{remark}
\textbf{Proof of Theorem~\ref{thm:freesubgroup}.} We may assume w.\,l.\,o.\,g.~that $\sphericalangle_{\{1,2\}}\geq\sphericalangle_{\{1,3\}}\geq\sphericalangle_{\{2,3\}}$. So, there are exactly three possibilities for the Gersten-Stallings angles, each of which is considered in a separate column in Figure~\ref{fig:ninecases}.

We distinguish between two cases. First, if there is an $a\in\{1,2,3\}$ such that the index of $\mathcal{G}|_{\emptyset}$ in $\mathcal{G}|_{\{a\}}$ is at least 3, consider the element $h\in\mathcal{G}$ that is given in the respective entry in Figure~\ref{fig:ninecases}. It is constructed in such a way that both $h$ and $h^{-1}$ are adapted to a billiard sequence $\mathcal{B}_{1}$ with the following property: The billiard shot $\pathbrackets{\mathcal{B}_{1}}$ starts at some point in the interior of $\Delta$ and goes orthogonally away from the edge labelled by $\{a\}$. After a couple of reflections, it comes back to the starting point, but in the opposite direction, see \circled{1} in Figure~\ref{fig:ninecases}.

Notice that, given an element $g\in\mathcal{G}$ with a decomposition into factors that are alternately from $\{h,h^{-1}\}$ and $\mathcal{G}|_{\{a\}}\setminus\mathcal{G}|_{\emptyset}$, we may concatenate the billiard shot $\pathbrackets{\mathcal{B}_{1}}$ and the orthogonal reflection at the edge labelled by $\{a\}$, see \circled{2} in Figure~\ref{fig:ninecases}, accordingly. This yields a billiard sequence $\mathcal{B}_{2}$, which $g$ is adapted to. Hence, we know: \emph{If there is at least one factor in the decomposition of $g$, then there is at least one reflection in the billiard sequence $\mathcal{B}_{2}$ and, by the billiards theorem, $g$ is non-trivial.}

Since the index of $\mathcal{G}|_{\emptyset}$ in $\mathcal{G}|_{\{a\}}$ is at least 3, we can find an element \mbox{$\widetilde{g}_{a}\in\mathcal{G}|_{\{a\}}$} that is neither in $\mathcal{G}|_{\emptyset}$ nor in $g_{a}\mathcal{G}|_{\emptyset}$. In particular, neither $\widetilde{g}_{a}{}^{-1}g_{a}$ nor ${g_{a}}^{-1}\widetilde{g}_{a}$ is in $\mathcal{G}|_{\emptyset}$. Define $x:=g_{a}h\widetilde{g}_{a}{}^{-1}\in\mathcal{G}$ and $y:=hg_{a}h\widetilde{g}_{a}{}^{-1}h^{-1}\in\mathcal{G}$. We claim that $x$ and $y$ generate a non-abelian free subgroup of $\mathcal{G}$.

Consider a non-empty freely reduced word over the letters $x$ and $y$ and their formal inverses. The element $g\in\mathcal{G}$ that is represented by this word has a natural decomposition into factors from $\{h^{\pm 1},{g_{a}}^{\pm 1},\widetilde{g}_{a}{}^{\pm 1}\}$. Cancel each $h^{-1}h$ and subsume each $\widetilde{g}_{a}{}^{-1}g_{a}$ and each ${g_{a}}^{-1}\widetilde{g}_{a}$ to a single element in $\mathcal{G}|_{\{a\}}\setminus\mathcal{G}|_{\emptyset}$. This yields a new decomposition of $g$ into factors that are alternately from $\{h,h^{-1}\}$ and $\mathcal{G}|_{\{a\}}\setminus\mathcal{G}|_{\emptyset}$. It is easy to see that, despite of the cancellation of each $h^{-1}h$, there is at least one factor left in the new decomposition of $g$. So, by our preliminary discussion, $g$ is non-trivial, which completes the proof that $x$ and $y$ generate a non-abelian free subgroup of $\mathcal{G}$.

Second, consider the case that there are two distinct $a,b\in\{1,2,3\}$ such that $\mathcal{G}|_{\{a,b\}}$ is not generated by $\mathcal{G}|_{\{a\}}$ and $\mathcal{G}|_{\{b\}}$. In this case, let $X:=\mathfrak{G}|_{\{a,b\}}$, $A:=\left\langle\mathfrak{G}|_{\{a\}},\mathfrak{G}|_{\{b\}}\right\rangle\leq\mathfrak{G}$, and $Y:=\left\langle\mathfrak{G}|_{\{a,c\}},\mathfrak{G}|_{\{b,c\}}\right\rangle\leq\mathfrak{G}$, where $c$ is the remaining element of $\{1,2,3\}$. Using the presentation (\ref{def:colimit}) one can show that $\mathcal{G}\cong X\ast_{A}Y$. Here, the homomorphisms are the ones induced by the inclusions.

By assumption, $|X:A|\geq 2$. On the other hand, we know that there is an element $g_{c}\in\mathcal{G}|_{\{c\}}\setminus\mathcal{G}|_{\emptyset}$. By the intersection theorem, $\mathcal{G}|_{\{c\}}\cap\mathcal{G}|_{\{a,b\}}=\mathcal{G}|_{\emptyset}$. So, $g_{c}\in\mathcal{G}|_{\{c\}}\setminus\mathcal{G}|_{\{a,b\}}\subseteq Y\setminus A$ and $|Y:A|\geq 2$. If $|Y:A|=2$, then $A$ is a normal subgroup of $Y$. Again, by the intersection theorem:
\[
\begin{array}{c}
{g_{c}}^{-1}g_{a}g_{c}\in A\cap\mathcal{G}|_{\{a,c\}}\subseteq\mathcal{G}|_{\{a,b\}}\cap\mathcal{G}|_{\{a,c\}}=\mathcal{G}|_{\{a\}} \medskip \\
{g_{c}}^{-1}g_{b}g_{c}\in A\cap\mathcal{G}|_{\{b,c\}}\subseteq\mathcal{G}|_{\{a,b\}}\cap\mathcal{G}|_{\{a,b\}}=\mathcal{G}|_{\{b\}}
\end{array}
\]
This implies that $\sphericalangle_{\{a,c\}}=\sphericalangle_{\{b,c\}}=\nicefrac{\pi\,}{2}$. Hence $\sphericalangle_{\{a,b\}}=0$, which is not possible since we assume the triangle of groups to be non-degenerate. So, $|Y:A|\geq 3$, and $\mathcal{G}\cong X\ast_{A}Y$ has a non-abelian free subgroup, see Remark~\ref{rem:bassserre}. \hfill$\Box$
\subsection{Tits alternative} \label{sub:titsalternative}
By Theorems~\ref{thm:hk} and \ref{thm:freesubgroup}, the colimit groups of many non-spherical triangles have a non-abelian free subgroup. In this section, we ask about the other cases and discuss the following version of the Tits alternative.
\begin{definition}[``Tits alternative'']
A class $\mathfrak{C}$ of groups satisfies the Tits alternative if each $G\in\mathfrak{C}$ is either large, i.\,e.~has a non-abelian free subgroup, or small, i.\,e.~is virtually solvable.
\end{definition}
\begin{remark}
There are groups that are neither large nor small. For example, take Thompson's group $F$. It has been shown by Brin and Squier in \cite{BrinSquier1985} that $F\leq\operatorname{PLF}(\mathbb{R})$ doesn't have a non-abelian free subgroup. And if $F$ was virtually solvable, then $[F,F]$ would have to be virtually solvable, too. But $[F,F]$ is infinite and simple, see \cite[Section 4]{CannonFloydParry1996}, which implies that $[F,F]$ cannot be virtually solvable.
\end{remark}
We may use Thompson's group $F$ to prove that the Tits alternative doesn't hold for the class of colimit groups of non-spherical triangles of groups. For further details about the following examples we refer to \cite{Cuno2011}.

Let $\Gamma_{1}$ be the triangle of groups with the property that the groups $G_{J}$ are all equal to $F$ and the injective homomorphisms $\varphi_{J_{1}J_{2}}:G_{J_{1}}\rightarrow G_{J_{2}}$ are all identities. The Gersten-Stallings angles amount to \mbox{$\sphericalangle_{\{1,2\}}=\sphericalangle_{\{1,3\}}=\sphericalangle_{\{2,3\}}=0$}, whence $\Gamma_{1}$ is a degenerate hyperbolic triangle of groups. But the colimit group is isomorphic to $F$ and, therefore, neither large nor small.

Notice that there are non-degenerate examples, too. Pick one of the three triangles of groups given in the introduction to Section~\ref{sec:billiards} and replace every group $G_{J}$ by $F\times G_{J}$ and every injective homomorphism $\varphi_{J_{1}J_{2}}:G_{J_{1}}\rightarrow G_{J_{2}}$ by $\operatorname{id}_{F}\times\varphi_{J_{1}J_{2}}:F\times G_{J_{1}}\rightarrow F\times G_{J_{2}}$. The Gersten-Stallings angles remain the same, whence the new triangle of groups, call it $\Gamma_{2}$, is non-degenerate and Euclidean. But the colimit group is isomorphic to $F\times\Delta(k,l,m)$ and, therefore, neither large nor small.\footnote{Recall the following two results, each of which is elementary: Let $G$ be a group and let $N\trianglelefteq G$ be a normal subgroup. $G$ has a non abelian free subgroup if and only if $N$ or $G\,/\,N$ does. So, $F\times\Delta(k,l,m)$ is not large. On the other hand, if $G$ is virtually solvable, then every subgroup of $G$ is virtually solvable, too. So, $F\times\Delta(k,l,m)$ is not small.}

Let us now assume that $G_{\emptyset}$ is finitely generated and either large or small, e.\,g.~$G_{\emptyset}=\{1\}$ as in Remark~\ref{rem:triviallyfilled}! In the non-degenerate case, this assumption already implies the Tits alternative.
\begin{theorem}
\label{thm:titsone}
The Tits alternative holds for the class of colimit groups $\mathcal{G}$ of non-degenerate non-spherical triangles of groups with the property that the group $G_{\emptyset}$ is finitely generated and either large or small.
\end{theorem}
Interestingly, in the degenerate case, it doesn't. Just consider the triangle of groups $\Gamma_{3}$ given by the following data:
\[
\begin{array}{c} G_{\varnothing}\cong\{1\},\;G_{\{1\}}\cong F,\;G_{\{2\}}\cong\left\langle a\,:\,a^{2}=1\right\rangle,\;G_{\{3\}}\cong\left\langle b\,:\,b^{2}=1\right\rangle, \smallskip \\ G_{\{1,2\}}\cong F\times\left\langle a\,:\,a^{2}=1\right\rangle,\;G_{\{1,3\}}\cong F\times\left\langle b\,:\,b^{2}=1\right\rangle,\;G_{\{2,3\}}\cong\left\langle a,b\,:\,a^{2}=b^{2}=1\right\rangle \end{array}
\]
Here, the homomorphisms $\varphi_{J_{1}J_{2}}:G_{J_{1}}\rightarrow G_{J_{2}}$ are given by $\forall\,f\in F:f\mapsto (f,1)$, by $a\mapsto (1,a)$ and $a\mapsto a$, and by $b\mapsto (1,b)$ and $b\mapsto b$. The Gersten-Stallings angles amount to $\sphericalangle_{\{1,2\}}=\nicefrac{\pi\,}{2}$, $\sphericalangle_{\{1,3\}}=\nicefrac{\pi\,}{2}$, $\sphericalangle_{\{2,3\}}=0$, whence $\Gamma_{3}$ is a degenerate Euclidean triangle of groups. Its colimit group is isomorphic to $F\times(\mathbb{Z}_{2}\ast\mathbb{Z}_{2})$ and, therefore, neither large nor small. The following theorem is an analogue of Theorem~\ref{thm:titsone}. It includes the degenerate case.
\begin{theorem}
\label{thm:titstwo}
The Tits alternative holds for the class of colimit groups $\mathcal{G}$ of non-spherical triangles of groups with the property that every group $G_{J}$ with $J\subseteq\{1,2,3\}$ and $|J|\leq 2$ is finitely generated and either large or small.
\end{theorem}
Given Theorem~\ref{thm:freesubgroup}, the proofs of Theorems~\ref{thm:titsone} and \ref{thm:titstwo} are elementary. But let us first discuss two auxiliary results!
\begin{lemma}
\label{lem:interchange}
Let $A$, $B$, and $C$ be groups such that $A$ contains $B$ as a normal subgroup with cyclic quotient and $B$ contains $C$ as a subgroup of finite index. If $B$ is finitely generated, then there is a subgroup $\widetilde{C}\leq C$ with the following property: There is a group $\widetilde{B}$ such that $A$ contains $\widetilde{B}$ as a subgroup of finite index and $\widetilde{B}$ contains $\widetilde{C}$ as a normal subgroup with cyclic quotient. In symbols:
\begin{center}
\input{img-21.pspdftex}
\end{center}
\end{lemma}
\textbf{Proof.} Pick an element $a\in A$ such that the cyclic quotient $A\,/\,B$ is generated by $aB$. Let $T$ be the respective \emph{transversal}, i.\,e.~the subgroup of $A$ that is generated by $a$.

It is well known, see e.\,g.~Marshall Hall's theorem \cite[Theorem 5.2]{Hall1949}, that the number of subgroups of index $n$ in a free group of finite rank $r$ is finite. So, by the correspondence theorem, see e.\,g.~\cite[Theorem 2.28]{Rotman1995}, the number of subgroups of index $|B:C|$ in the finitely generated group $B$ is finite, too. Let $\widetilde{C}$ be the intersection of all these subgroups, whence $\widetilde{C}\leq C\leq B$. Moreover, $\widetilde{C}$ is a subgroup of finite index in $B$, see e.\,g.~\cite[Exercise 3.31\,(i)]{Rotman1995}, and it is \emph{characteristic} in $B$, i.\,e.~it is fixed by any automorphism $\varphi:B\rightarrow B$. For more details, see also \cite[Satz 1.15]{CampsRebelRosenberger2008}.

Since $A$ contains $B$ as a normal subgroup, conjugation by any element from $A$ induces an automorphism $\varphi:B\rightarrow B$. But $\widetilde{C}$ is characteristic in $B$. So, $\widetilde{C}$ is fixed by $\varphi$, whence $A$ contains $\widetilde{C}$ as a normal subgroup. Therefore, the set $\widetilde{B}:=T\widetilde{C}=\{\,t\widetilde{c}\mid t\in T,\widetilde{c}\in\widetilde{C}\,\}\subseteq A$ is actually a subgroup. Moreover, $\widetilde{B}=T\widetilde{C}$ contains $\widetilde{C}$ as a normal subgroup, and the quotient $\widetilde{B}\,/\,\widetilde{C}=T\widetilde{C}\,/\,\widetilde{C}\cong T\,/\,(T\cap\widetilde{C})$ is cyclic.

We still need to show that $A$ contains $\widetilde{B}$ as a subgroup of finite index. Recall that for every $g\in A$ there is an exponent $k\in\mathbb{Z}$ such that $g\in a^{k}B$. Since $\widetilde{C}$ is a subgroup of finite index in $B$, we can find a finite set $\{b_{1},b_{2},\dotsc,b_{m}\}$ of representatives of the right cosets of $\widetilde{C}$ in $B$. So, there is an $i\in\{1,2,\dotsc,m\}$ such that $g\in a^{k}\widetilde{C}b_{i}\subseteq T\widetilde{C}b_{i}=\widetilde{B}b_{i}$. In other words, every $g\in A$ is contained in one of the cosets $\widetilde{B}b_{1},\widetilde{B}b_{2},\dotsc,\widetilde{B}b_{m}$, which finishes the proof. \hfill$\Box$%
\begin{corollary}
\label{cor:virtuallysolvable}
Let $G$ be a group and let $H$ be a normal subgroup of $G$. If $G\,/\,H$ is virtually polycyclic and $H$ is finitely generated and virtually solvable, then $G$ is virtually solvable.
\end{corollary}
\textbf{Proof.} Since $G\,/\,H$ is virtually polycyclic, there is a series:
\begin{center}
\input{img-22.pspdftex}
\end{center}
By the correspondence theorem, see e.\,g.~\cite[Theorem 2.28]{Rotman1995}, this yields a series from $G$ to $H$ with the same, i.\,e.~isomorphic, quotients:
\begin{center}
\input{img-23.pspdftex}
\end{center}
We may concatenate the latter with the series from $H$ to $\{1\}$ and obtain:
\begin{center}
\input{img-24.pspdftex}
\end{center}
Since $H$ is finitely generated, we may now apply Lemma~\ref{lem:interchange} and replace:
\begin{center}
\input{img-25.pspdftex}
\end{center}
Finally, we intersect the lower terms $C_{2},C_{3},\dotsc,C_{n}$ with $\widetilde{C}_{1}$. This gives us a new series from $G$ to $\{1\}$ with the same structure: a finite index subgroup, cyclic quotients, a finite index subgroup, abelian quotients. Of course, the position of the second finite index subgroup has moved one step to the left.

Since the cyclic extensions $B_{m},B_{m-1},\dotsc,B_{1}$ are also finitely generated, we may iterate the above procedure until the two finite index subgroups are subsequent, which proves that $G$ is virtually solvable. \hfill $\Box$\bigskip\\
\textbf{Proof of Theorem~\ref{thm:titsone}.} Consider a non-degenerate non-spherical triangle of groups with the property that the group $G_{\emptyset}$ is finitely generated and either large or small. If $G_{\emptyset}$, and hence $\mathcal{G}|_{\emptyset}$, is large, then the colimit group $\mathcal{G}$ is large, too. So, we may assume w.\,l.\,o.\,g.~that $G_{\emptyset}$, and hence $\mathcal{G}|_{\emptyset}$, is small.

The triangle of groups is non-degenerate. So, for every $a\in\{1,2,3\}$ there is an element $g_{a}\in\mathcal{G}|_{\{a\}}\setminus\mathcal{G}|_{\emptyset}$. If the triangle of groups is hyperbolic, then, by Theorem~\ref{thm:hk}, the colimit group $\mathcal{G}$ is large. So, we may assume w.\,l.\,o.\,g.~that the triangle of groups is Euclidean. Now, let $a\in\{1,2,3\}$. Since $g_{a}\in\mathcal{G}|_{\{a\}}\setminus\mathcal{G}|_{\emptyset}$, the index of $\mathcal{G}|_{\emptyset}$ in $\mathcal{G}|_{\{a\}}$ is at least $2$. If it is strictly larger than $2$, then the simplicial complex $\mathcal{X}$ branches and, by Theorem~\ref{thm:freesubgroup}, $\mathcal{G}$ is large. So, we may assume w.\,l.\,o.\,g.~that it is equal to $2$. In particular, $\mathcal{G}|_{\emptyset}$ is normal in $\mathcal{G}|_{\{a\}}$. By the same argument, we may assume w.\,l.\,o.\,g.~that for every two distinct $a,b\in\{1,2,3\}$ the group $\mathcal{G}|_{\{a,b\}}$ is generated by $\mathcal{G}|_{\{a\}}$ and $\mathcal{G}|_{\{b\}}$. Therefore, $\mathcal{G}|_{\emptyset}$ is normal in $\mathcal{G}|_{\{1\}}$, $\mathcal{G}|_{\{2\}}$, $\mathcal{G}|_{\{3\}}$, $\mathcal{G}|_{\{1,2\}}$, $\mathcal{G}|_{\{1,3\}}$, $\mathcal{G}|_{\{2,3\}}$, and, finally, in $\mathcal{G}$.

Notice that this property does also hold in the triangle of groups itself, i.\,e.~for every $J\subseteq\{1,2,3\}$ with $1\leq |J|\leq 2$ the image $\varphi_{\emptyset J}(G_{\emptyset})$ is normal in $G_{J}$. For a formal proof, apply the natural homomorphism $\nu_{J}:G_{J}\rightarrow\mathcal{G}$, which is injective, and observe:
\[
\begin{array}{r@{\;}c@{\;}l}
\nu_{J}\circ\varphi_{\emptyset J}(G_{\emptyset}) & = & \nu_{\emptyset}(G_{\emptyset}) = \tilde\nu_{\emptyset}\circ\mu_{\emptyset}(G_{\emptyset}) = \tilde\nu_{\emptyset}(\mathcal{G}|_{\emptyset}) \smallskip \\ & \trianglelefteq & \tilde\nu_{J}(\mathcal{G}|_{J}) = \tilde\nu_{J}\circ\mu_{J}(G_{J}) = \nu_{J}(G_{J}) \\
\end{array}
\]
We may therefore construct the \emph{quotient triangle of groups}, which is obtained by replacing the group $G_{\emptyset}$ by $G_{\emptyset}\,/\,G_{\emptyset}\cong\{1\}$ and for every $J\subseteq\{1,2,3\}$ with $1\leq |J|\leq 2$ the group $G_{J}$ by its quotient $G_{J}\,/\,\varphi_{\emptyset J}(G_{\emptyset})$.

At this point, one needs to verify that every injective homomorphism $\varphi_{J_{1}J_{2}}:G_{J_{1}}\rightarrow G_{J_{2}}$ induces an injective homomorphism between the quotients and that the Gersten-Stallings angles remain the same. We leave this work to the reader. However, the resulting diagram is a non-degenerate Euclidean triangle of groups. Moreover, using the presentation $(\ast)$ one can show that its colimit group is isomorphic to $\mathcal{G}\,/\,\mathcal{G}|_{\emptyset}$. We will now study the quotient triangle of groups in some more detail.

For every $a\in\{1,2,3\}$ the index of $\varphi_{\emptyset\{a\}}(G_{\emptyset})$ in $G_{\{a\}}$ is equal to $2$, i.\,e.~the quotient $G_{\{a\}}\,/\,\varphi_{\emptyset\{a\}}(G_{\emptyset})$ has two elements. Therefore, the quotient triangle of groups must be isomorphic to one of the three triangles of groups given in the introduction to Section~\ref{sec:billiards} and, in particular, its colimit group $\mathcal{G}\,/\,\mathcal{G}|_{\emptyset}$ must be isomorphic to a subgroup $S\leq\operatorname{Isom}(\mathbb{E}^{2})$. Let $\operatorname{Trans}(\mathbb{E}^{2})\leq\operatorname{Isom}(\mathbb{E}^{2})$ be the subgroup of all translations. We claim that the intersection $T:=S\cap\operatorname{Trans}(\mathbb{E}^{2})$ is of finite index in $S$ and isomorphic to $\mathbb{Z}^{2}$ or $\mathbb{Z}$ or $\{1\}$, in fact, it will be $\mathbb{Z}^{2}$. This means that $S$ is virtually polycyclic. And so is the colimit group $\mathcal{G}\,/\,\mathcal{G}|_{\emptyset}$ of the quotient triangle of groups.

Let us prove the claim! First, let $\operatorname{Isom}^{+}(\mathbb{E}^{2})\leq\operatorname{Isom}(\mathbb{E}^{2})$ be the subgroup of all orientation preserving isometries and observe that $S\cap\operatorname{Isom}^{+}(\mathbb{E}^{2})\leq S$ is of index $2$. Second, observe that $T\leq S\cap\operatorname{Isom}^{+}(\mathbb{E}^{2})$ is of finite index, too. Indeed, every non-trivial coset of $T\leq S\cap\operatorname{Isom}^{+}(\mathbb{E}^{2})$ consists of all rotations in $S$ with the same angle of rotation. But it is easy to see that every possible angle of rotation can be expressed as $\lambda_{1}\cdot\nicefrac{2\pi\,}{k}+\lambda_{2}\cdot\nicefrac{2\pi\,}{l}+\lambda_{3}\cdot\nicefrac{2\pi\,}{m}$ with $\lambda_{i}\in\mathbb{Z}$, whence their total number is finite. Third, by construction of $S$, the images of the triangle enclosed by the three lines yield a tessellation of $\mathbb{E}^{2}$. So, the orbit of each point $p\in\mathbb{E}^{2}$ is a discrete subset of $\mathbb{E}^{2}$. In particular, $T\leq\operatorname{Trans}(\mathbb{E}^{2})\cong\mathbb{R}^{2}$ is a discrete subgroup and, hence, isomorphic to $\mathbb{Z}^{2}$ or $\mathbb{Z}$ or $\{1\}$.

So, we know that the colimit group $\mathcal{G}\,/\,\mathcal{G}|_{\emptyset}$ of the quotient triangle of groups is virtually polycyclic. Moreover, by assumption, $G_{\varnothing}$, and hence $\mathcal{G}|_{\varnothing}$, is finitely generated and virtually solvable. So, by Corollary~\ref{cor:virtuallysolvable}, we may conclude that $\mathcal{G}$ is virtually solvable, i.\,e.~small. \hfill$\Box$\bigskip\\
\textbf{Proof of Theorem 10.} Consider a non-spherical triangle of groups with the property that every group $G_{J}$ with $J\subseteq\{1,2,3\}$ and $|J|\leq 2$ is finitely generated and either large or small. Again, we may assume w.\,l.\,o.\,g.~that the groups $G_{J}$, and hence their images $\mathcal{G}|_{J}$, are small.

Moreover, if the triangle of groups is non-degenerate, then we know by Theorem~\ref{thm:titsone} that $\mathcal{G}$ is either large or small. So, we may assume w.\,l.\,o.\,g.~that the Gersten-Stallings angle $\sphericalangle_{\{2,3\}}=0$, which means that the homomorphism $\alpha:G_{\{2\}}\ast_{G_{\varnothing}}G_{\{3\}}\rightarrow G_{\{2,3\}}$ induced by $\varphi_{\{2\}\{2,3\}}$ and $\varphi_{\{3\}\{2,3\}}$ is injective.

As already mentioned in the proof of Theorem~\ref{thm:freesubgroup}, one can always show that $\mathcal{G}\cong X\ast_{A}Y$ with $X:=\mathcal{G}|_{\{2,3\}}$, $A:=\left\langle\mathcal{G}|_{\{2\}},\mathcal{G}|_{\{3\}}\right\rangle\leq\mathcal{G}$, $Y:=\left\langle\mathcal{G}|_{\{1,2\}},\mathcal{G}|_{\{1,3\}}\right\rangle\leq\mathcal{G}$.  Depending on $|X:A|$ and $|Y:A|$, we distinguish between four cases:

\begin{enumerate}
\item If $|X:A|=1$, then $\mathcal{G}|_{\{2,3\}}$ is generated by $\mathcal{G}|_{\{2\}}$ and $\mathcal{G}|_{\{3\}}$ or, equivalently, $G_{\{2,3\}}$ is generated by $\varphi_{\{2\}\{2,3\}}(G_{\{2\}})$ and $\varphi_{\{3\}\{2,3\}}(G_{\{3\}})$. So, $\alpha$ is not only injective but also surjective, whence $G_{\{2\}}\ast_{G_{\varnothing}}G_{\{3\}}\cong G_{\{2,3\}}$. This allows us to simplify the original presentation (\ref{def:colimit}) of the colimit group $\mathcal{G}$ by deleting superficial generators and relators:
\[
\begin{array}{r@{~}l}
\mathcal{G}\cong\big\langle G_{\{1\}},G_{\{1,2\}},G_{\{1,3\}}\,: & R_{\{1\}},R_{\{1,2\}},R_{\{1,3\}}, \smallskip \\
& \left\{g=\varphi_{\{1\}\{1,2\}}(g)\,:\,g\in G_{\{1\}}\right\}, \smallskip \\
& \left\{g=\varphi_{\{1\}\{1,3\}}(g)\,:\,g\in G_{\{1\}}\right\}\big\rangle
\end{array}
\]
So, $\mathcal{G}\cong G_{\{1,2\}}\ast_{G_{\{1\}}}G_{\{1,3\}}$ or, equivalently, $\mathcal{G}\cong\mathcal{G}|_{\{1,2\}}\ast_{\mathcal{G}|_{\{1\}}}\mathcal{G}|_{\{1,3\}}$. Now, we may, again, distinguish between four cases: 
\begin{enumerate}
\item If $|\,\mathcal{G}|_{\{1,2\}}:\mathcal{G}|_{\{1\}}\,|=1$, then $\mathcal{G}\cong\mathcal{G}|_{\{1,3\}}$, which is small.
\item If $|\,\mathcal{G}|_{\{1,3\}}:\mathcal{G}|_{\{1\}}\,|=1$, then $\mathcal{G}\cong\mathcal{G}|_{\{1,2\}}$, which is small.
\item If $|\,\mathcal{G}|_{\{1,2\}}:\mathcal{G}|_{\{1\}}\,|=|\,\mathcal{G}|_{\{1,3\}}:\mathcal{G}|_{\{1\}}\,|=2$, then $\mathcal{G}|_{\{1\}}$ is normal in $\mathcal{G}|_{\{1,2\}}$, $\mathcal{G}|_{\{1,3\}}$, and $\mathcal{G}$. The quotient $\mathcal{G}\,/\,\mathcal{G}|_{\{1\}}\cong\mathbb{Z}_{2}\ast\mathbb{Z}_{2}$, which is virtually (poly-)cyclic. On the other hand, $\mathcal{G}|_{\{1\}}$ itself is finitely generated and virtually solvable. So, by Corollary~\ref{cor:virtuallysolvable}, the colimit group $\mathcal{G}$ is virtually solvable, i.\,e.~is small.
\item Otherwise, $|\,\mathcal{G}|_{\{1,2\}}:\mathcal{G}|_{\{1\}}\,|\geq 2$ and $|\,\mathcal{G}|_{\{1,3\}}:\mathcal{G}|_{\{1\}}\,|\geq 2$ and not both equal to 2. But then, by Remark~\ref{rem:bassserre}, $\mathcal{G}\cong\mathcal{G}|_{\{1,2\}}\ast_{\mathcal{G}|_{\{1\}}}\mathcal{G}|_{\{1,3\}}$ has a non-abelian free subgroup, i.\,e.~is large.
\end{enumerate}
\item If $|Y:A|=1$, then $\mathcal{G}\cong X=\mathcal{G}|_{\{2,3\}}$, which is small.
\item If $|X:A|=|Y:A|=2$, then $A$ is normal in $X$, $Y$, and $\mathcal{G}$. The quotient $\mathcal{G}\,/\,A\cong\mathbb{Z}_{2}\ast\mathbb{Z}_{2}$, which is virtually (poly-)cyclic. Let us now study $A$ in more detail! Since $\alpha:G_{\{2\}}\ast_{G_{\varnothing}}G_{\{3\}}\rightarrow G_{\{2,3\}}$ is injective, the subgroups $\left\langle\varphi_{\{2\}\{2,3\}}(G_{2}),\varphi_{\{3\}\{2,3\}}(G_{3})\right\rangle\leq G_{\{2,3\}}$ and $A=\left\langle\mathcal{G}|_{\{2\}},\mathcal{G}|_{\{3\}}\right\rangle\leq\mathcal{G}|_{\{2,3\}}\leq\mathcal{G}$ are both isomorphic to $G_{\{2\}}\ast_{G_{\varnothing}}G_{\{3\}}$. In particular, $A$ is finitely generated. By a case analysis analogue to the one in (1), we can see that either $A$ is virtually solvable, namely in cases (a)--(c), or $A$ has a non-abelian free subgroup, namely in case (d). In the former cases, by Corollary~\ref{cor:virtuallysolvable}, the colimit group $\mathcal{G}$ is virtually solvable, i.\,e.~is small. In the latter case, the colimit group $\mathcal{G}$, which contains $A$ as a subgroup, has a non-abelian free subgroup, i.\,e.~is large.
\item Otherwise, $|X:A|\geq 2$ and $|Y:A|\geq 2$ and not both equal to 2. But then, by Remark~\ref{rem:bassserre}, $\mathcal{G}\cong X\ast_{A}Y$ has a non-abelian free subgroup, i.\,e.~is large. \hfill$\Box$
\end{enumerate}
\begin{remark}
Kopteva and Williams have proved the Tits alternative for the class of non-spherical Pride groups based on graphs with at least four vertices, see \cite[Theorem 1]{KoptevaWilliams2008}. One way to read Theorem~\ref{thm:titsone} is the following: The Tits alternative does not hold for the class of non-spherical Pride groups based on graphs with three vertices. But once we assume that each edge is genuine, i.\,e.~that none of the Gersten-Stallings angles is equal to 0, it does.
\end{remark}
%
% ==================================================
% Anhang
% ==================================================
%
\section*{Appendix}
In Sections~\ref{sub:bridson} and \ref{sub:billiardstheorem}, we have been working with non-degenerate Euclidean triangles of groups. But the construction can be extended to non-degenerate non-spherical triangles of groups. In the hyperbolic case, one can either pick a triangle $\Delta$ in the hyperbolic plane $\mathbb{H}^{2}$, as suggested by Bridson in \cite[p.\,431, ll.\,13--16]{Bridson1991}, or a triangle $\Delta$ in the Euclidean plane $\mathbb{E}^{2}$ ``whose angles are perhaps a little bit larger than the group-theoretic angles,'' as suggested by Gersten and Stallings in \cite[p.\,499, ll.\,7--9]{Stallings1991}. Let us sketch an application for each of the two alternatives!
\subsection*{A word about normal forms}
In the billiards theorem, we assume given an element $g\in\mathcal{G}$ that is adapted to a billiard sequence, i.\,e.~that is equipped with a suitable decomposition into factors from $\mathcal{G}|_{\{1\}}\setminus\mathcal{G}|_{\emptyset}$, $\mathcal{G}|_{\{2\}}\setminus\mathcal{G}|_{\emptyset}$, $\mathcal{G}|_{\{3\}}\setminus\mathcal{G}|_{\emptyset}$. But we could also go the other way and use the simplicial complex $\mathcal{X}$ to construct decompositions.

More precisely, given a non-degenerate non-spherical triangle of groups, pick a triangle $\Delta$ either in the Euclidean plane $\mathbb{E}^{2}$ or in the hyperbolic plane $\mathbb{H}^{2}$ whose angles agree with the Gersten-Stallings angles of the triangle of groups and construct the simplicial complex $\mathcal{X}$. Notice that all the results from Sections~\ref{sub:bridson} and \ref{sub:billiardstheorem} still hold true. Given an arbitrary element $g\in\mathcal{G}$, consider the unique geodesic, see \cite[Section 2, Main Theorem, (2)\,$\Rightarrow$\,(1)]{Bridson1991}, from the barycentre of $\sigma$ to the barycentre of $g\sigma$ in $|\mathcal{X}|$. As soon as $g\not\in\mathcal{G}|_{\emptyset}$, the geodesic traverses several $2$-simplices. First, it traverses $\sigma$. Then, depending on whether it leaves $\sigma$ crossing a $0$-simplex or the interior of a $1$-simplex, there is an element $g_{1}\in\mathcal{G}|_{K}$ with $K\subseteq\{1,2,3\}$ and $|K|=2$ or $|K|=1$ such that the next $2$-simplex it traverses is $g_{1}\sigma$.

\begin{figure}
\begin{center}
\input{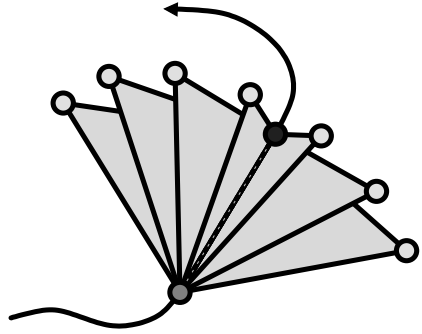}
\end{center}
\caption{There are many possibilities to choose a $2$-simplex.}
\label{fig:pages}
\end{figure}

This procedure goes on. At the end, it traverses $g_{1}g_{2}\dotsb g_{m}\sigma=g\sigma$, which yields a decomposition of $g$ into factors from the groups $\mathcal{G}|_{K}$ with $K\subseteq\{1,2,3\}$ and $1\leq |K|\leq 2$ and one final factor from $\mathcal{G}|_{\emptyset}$. Notice that this decomposition is not well defined, even in the case $\mathcal{G}|_{\emptyset}=\{1\}$. Just imagine the geodesic running along some $1$-simplex. Then, there are many possibilities to choose the respective $2$-simplex $g_{1}g_{2}\dotsb g_{i}\sigma$, see Figure~\ref{fig:pages}. On the other hand, if we fix a set of coset representatives for each pair of subgroups $\mathcal{G}|_{K_{1}}\leq\mathcal{G}|_{K_{2}}$ with $K_{1}\subset K_{2}\subseteq\{1,2,3\}$ and $|K_{2}|\leq 2$, there is a well defined decomposition in terms of these coset representatives and one final factor from $\mathcal{G}|_{\emptyset}$. Even though it seems to be inconvenient to work with, we may call it a normal form.
\subsection*{Euclidean domination}
The second alternative has the advantage that we need to study only three different triangles $\Delta$. More precisely, given a non-degenerate non-spherical triangle of groups, the Gersten-Stallings angles are always of the form $\nicefrac{2\pi\,}{\hat{m}}$, where $\hat{m}$ is even. Let us think of them as $\nicefrac{\pi\,}{k}$, $\nicefrac{\pi\,}{l}$, $\nicefrac{\pi\,}{m}$ with $k,l,m\in\mathbb{N}$ and $k\leq l\leq m$!
\begin{lemma} \label{lem:Euclideandomination}
Let $k,l,m\in\mathbb{N}$ as above, i.\,e.~with $k\leq l\leq m$ and $\nicefrac{\pi\,}{k}+\nicefrac{\pi\,}{l}+\nicefrac{\pi\,}{m}\leq\pi$. Then, there are $k',l',m'\in\mathbb{N}$ with $k'\leq k$, $l'\leq l$, $m'\leq m$ and $\nicefrac{\pi\,}{k'}+\nicefrac{\pi\,}{l'}+\nicefrac{\pi\,}{m'}=\pi$.
\end{lemma}
\textbf{Proof.} By assumption, $k\geq 2$. If $k\geq 3$, let $k'=l'=m'=3$. Now, consider the case $k=2$. Then, by assumption, $l\geq 3$. If $l\geq 4$, let $k'=2$ and $l'=m'=4$. On the other hand, if $l=3$, then $m\geq 6$ and we may choose $k'=2$, $l'=3$, $m'=6$. \hfill$\Box$\bigskip \\
So, either $(\nicefrac{\pi\,}{3},\nicefrac{\pi\,}{3},\nicefrac{\pi\,}{3})$ or $(\nicefrac{\pi\,}{2},\nicefrac{\pi\,}{4},\nicefrac{\pi\,}{4})$ or $(\nicefrac{\pi\,}{2},\nicefrac{\pi\,}{3},\nicefrac{\pi\,}{6})$ \emph{dominates} $(\nicefrac{\pi\,}{k},\nicefrac{\pi\,}{l},\nicefrac{\pi\,}{m})$, i.\,e.~is  coordinatewise at least as large as $(\nicefrac{\pi\,}{k},\nicefrac{\pi\,}{l},\nicefrac{\pi\,}{m})$. If we take this dominating triple instead of the original Gersten-Stallings angles, then, again, all results from Sections~\ref{sub:bridson} and \ref{sub:billiardstheorem}, in particular the link condition and the billiards theorem, hold true. This way, one can easily see that the proof of Theorem~\ref{thm:freesubgroup} extends to all non-degenerate non-spherical triangles of groups.

% ==================================================
% BIBLIOGRAPHIE
% ==================================================
%

%
\end{document}